\documentclass{article}

\usepackage{arxiv}

\usepackage{amsfonts}       
\usepackage{xcolor}

\usepackage[authoryear]{natbib}
\usepackage{hyperref}

\usepackage{paralist}

\usepackage{setspace}
\usepackage{amsmath}
\usepackage{amssymb}
\usepackage{mathtools}
\usepackage{mismath}
\usepackage{derivative}
\usepackage{subcaption}
\usepackage{lineno}

\usepackage{float}
\usepackage{booktabs}
\usepackage{paralist}
\usepackage{array}
\usepackage{algpseudocode, algorithm, algorithmicx}

\DeclarePairedDelimiter{\pthTemp}{(}{)}
\DeclarePairedDelimiter{\brkTemp}{[}{]}
\DeclarePairedDelimiter{\bcolTemp}{\{}{\}}

\newcommand{\say}[1]{``{#1}''}
\newcommand{\pth}[1]{\pthTemp{#1}}
\newcommand{\brk}[1]{\brkTemp{#1}}
\newcommand{\bcol}[1]{\bcolTemp{#1}}

\newcommand{\dbcol}[1]{\bcol{\!\! \bcol{#1} \!\! }}
\newcommand{\dbrk}[1]{\brk{\!\brk{#1}\!}}
\DeclarePairedDelimiter{\seminrm}{\lvert}{\rvert}
\newcommand{\nrm}[1]{\seminrm{ \hspace{-0.4mm}\seminrm{#1} \hspace{-0.4mm}}}

\newcommand{\ra}[1]{\renewcommand{\arraystretch}{#1}}

\newcommand{\dx}{\,{\rm d} \mathbf{x}}
\newcommand{\dS}{\,{\rm d} S}
\newcommand{\meansc}[1]{\{\!\!\{{#1}\}\!\!\}}
\newcommand{\jumpsc}[1]{[\!\![{#1}]\!\!]}
\newcommand{\ndg}[1]{|\!|\!|{#1}|\!|\!|}

\title{Parallel Nodal Interior-Penalty Discontinuous Galerkin Methods for the
	Subsonic Compressible Navier-Stokes Equations: Applications to
	Vortical Flows and VIV Problems}

\author{
 Spiros Zafeiris \\
  School of Naval Architecture and Marine Engineering\\
  National Technical University of Athens\\
  Zografos, 15780, Athens \\
   \And
 Emmanuil H. Georgoulis \\
  School of Applied Mathematical and Physical Sciences\\
  National Technical University of Athens\\
  Zografos, 15780, Athens \\
  The Maxwell Institute for Mathematical Sciences \\ \& Department of Mathematics\\
  Heriot-Watt University \\
  Edinburgh, EH4 4AS, United Kingdom\\
  Institute for Applied and Computational Mathematics\\
  Foundation for Research and Technology -- Hellas\\
  Heraklion, Crete, Greece
  \And
 George Papadakis \\
  School of Naval Architecture and Marine Engineering\\
  National Technical University of Athens\\
  Zografos, 15780, Athens \\
}

\begin{document}
\boldvect
\maketitle
\begin{abstract}
    We present a Discontinuous Galerkin (DG) solver for the compressible
	Navier-Stokes system, designed for applications of technological and
	industrial interest in the subsonic region. More precisely, this work aims to exploit the
	DG-discretised Navier-Stokes for two dimensional vortex-induced vibration (VIV) problems allowing for high-order of accuracy.
	The numerical discretisation comprises a nodal DG method on triangular grids, that includes two
	types of numerical fluxes: 1) the Roe approximate Riemann solver flux for non-linear advection terms, and
	2) an Interior-Penalty numerical flux for non-linear diffusion terms.
	The nodal formulation permits the use of high order polynomial approximations without
	compromising computational robustness. The spatially-discrete form is integrated in time using a low-storage strong-stability preserving explicit Runge-Kutta scheme, and is coupled weakly with an implicit rigid body dynamics algorithm. The proposed algorithm successfully implements polynomial orders of $p\ge 4$ for the laminar compressible Navier-Stokes equations in massively parallel architectures. The resulting framework is firstly tested in terms of its convergence properties. Then, numerical solutions are validated with experimental and numerical data for the case of a circular cylinder at low Reynolds number, and lastly, the methodology is employed to simulate the problem of an  elastically-mounted cylinder, a known configuration characterised by significant computational challenges. The above results showcase that the DG framework can be employed as an accurate and efficient, arbitrary order numerical methodology, for high fidelity fluid-structure-interaction (FSI) problems.
\end{abstract}



\maketitle
\section{Introduction}

In the modern era of computational modelling, high-fidelity numerical methods have become essential for both designing engineering structures as well as understanding complex physical phenomena. 
Among these methods, the discontinuous Galerkin (DG) approach has emerged as particularly promising for fluid dynamics simulations because of its ability to achieve high-order accuracy on unstructured meshes. 
DG methods were first introduced by Reed and Hill \cite{reed1973triangular} for the neutron transport equation, and have since gained popularity for their effectiveness in modelling convection-dominated problems.
Error analysis for arbitrary order DG methods for the linear transport equation on unstructured meshes was established in \cite{johnson1986analysis}. Following, the development of Runge-Kutta DG (RKDG) methods in the late 1980s-1990s by Cockburn \& Shu and co-workers \cite{cockburn1989tvbII,cockburn1989tvbIII,cockburn1990tvbIV,cockburn1998tvbV} positioned prominently DG methodologies for computational fluid dynamics (CFD), and ever since, DG methods for the compressible Euler system have been developed in the 1990s and early 2000s \cite{BassiRebay1997Euler,vanderVegt2002,HartmannHouston2002,RemacleFlaherty,wang2007implicit}, followed by respective developments for the compressible Navier-Stokes system \cite{bassi1997high,hartmann2006adaptive,hartmann2008optimal,klaij2006space,gassner2008discontinuous} and later references.

A DG method can be locally variable order, can incorporate unstructured grids of simplicial, box-type, or even general polygonal/polyhedral \cite{CDGH14,CDGH17} element shapes, has local conservation properties and integrates seamlessly with stable numerical fluxes/approximate Riemann solvers and slope limiters. These features render the family of DG methods as a highly versatile numerical framework for studying complex flows with wake structures, whereby the preservation and accurate propagation of flow information downstream is critical. In such scenarios, typical and widely used second-order methods often introduce excessive numerical dissipation that artificially dampens vortical structures.

A particularly challenging example of flows with wake structures is the class of vortex-induced vibration (VIV) problems. A VIV problem is described by an immersed bluff body which is free to move under some elastic constraints. 
Such constraints are linear or rotational springs and dampers. In this coupled system, vortices shed by the solid body cause movement around the equilibrium. This motion then modifies the ongoing vortex shedding process and the result is a non-linear oscillator in which a \textit{lock-in} region is created (see \cite{SarpkayaVIO}). Numerical frameworks implemented on VIV problems include finite element \cite{mittal1992finite,mittal1999finite}, immersed boundary \cite{griffith2017sharp}, and hybrid Lagrangian-Eulerian \cite{papadakis2022hybrid} methods.

Somewhat surprisingly, there is extremely limited work on the development of DG methods for Vortex-Induced Vibration (VIV) problems of bluff bodies. This is despite the aforementioned pertinent properties of discontinuous Galerkin approaches in terms of extreme mesh flexibility, variable order incorporation and local conservation. Indeed, we are only aware of the very recent preprint \cite{zou2025moving} employing RKDG with interior penalty for a number of VIV benchmarks. 

This work is concerned with the development of a RKDG method for the time-dependent compressible Navier-Stokes system. The viscous part is discretised through an interior penalty (IP) flux \cite{baker1977finite,wheeler1978elliptic,arnold1982interior}, to take advantage of its compact stencil, thereby allowing for minimal communication on parallel implementations. Other established DG approaches include the Local DG method \cite{cockburn1998local} and the original methods of Bassi \& Rebay \cite{bassi1997high,bassi2005discontinuous}, which require wider stencils, or more recent approaches such as the compact DG \cite{PeraireCDG}, the hybrid DG \cite{CockburnHDG} and the embedded DG \cite{NguyenEHDG} methods. 
One the other hand, IP-DG requires a user-defined discontinuity penalization (penalty) parameter $\sigma=Cp^2/h$ with $h$ denoting the local mesh size, $p$ the local polynomial order, and $C>0$ typically selected ``large enough'' heuristically; this is also the choice of $\sigma$ in \cite{zou2025moving}. 
Large values of $\sigma$, however, reduce the admissible time-step size through the CFL condition; cf., \eqref{eq:cfl_general} below. 
Hence, the accurate definition of $\sigma$ is crucial to achieve efficient long-time simulations. 
To that end, a central development in this work is the accurate definition of \emph{sufficiently small} $\sigma$, without compromising the stability of the method. 
The new definition of $\sigma$, apart from $h$ and $p$ is also sensitive to the local direction magnitude of the non-linear diffusion. 
This choice of $\sigma$, in turn, allows for larger time-steps, thereby reducing the computation cost.  

In addition, the implementation takes advantage of the nodal basis framework from \cite{hesthaven2008nodal} employing optimal interpolation points \cite{bendito2007estimation,chen1995approximate,hesthaven1998electrostatics,taylor2000algorithm,gassner2009polymorphic} for the reduced conditioning of local mass matrices. The choice of nodal DG basis is motivated by the efficient calculation of the elemental neighbour contributions. The assembly step of the local stiffness matrices is performed through the Dubiner-Proriol-Koornwinder construction \cite{dubiner1991spectral,proriol1957famille,koornwinder1975two} for triangular elements. 
Affine elemental maps are used throughout, giving constant elemental Jacobian matrices and, thus, integrals are calculated and stored only once in the reference space for all elements. Of course, if desired, the proposed scheme can also accept curvilinear elements \cite{hesthaven2008nodal,winters2021construction}. Lastly, upwinding is implemented through Roe's approximate Riemann solver \cite{roe1981approximate} which is the standard for the compressible Navier-Stokes system, whereas other standard numerical fluxes \cite{toro2013riemann} are also admissible.

The resulting algorithm of the nodal IP RKDG method is implemented in parallel using the Message Passing Interface. For this reason, we take advantage of the stencil compactness offered by the IP flux and the nodal basis, and we report the resulting parallel efficiency. 
Also, the method is assessed for the solution of the subsonic compressible system \eqref{eq:cnp} on the basis of the convergence rate when using the newly defined interior penalty.

The remainder of this work is structured as follows. In the Section \ref{section:discretisation}, we present the numerical formulation of the method including the space and time strategies, the rigid body dynamics equation and the algorithm that couples them. The third section is focused on computational results, which are comprised of a manufactured solution, the von-K\'{a}rm\'{a}n vortex street of non-vibrating circular cylinder and the elastically-mounted cylinder. Finally, we draw some conclusions and summarise the findings.

\section{Discretisation}

\label{section:discretisation}
\newcommand{\ver}{\varphi}
\newcommand{\pvect}[1]{\pth{\vect{#1}}}
\newcommand{\vand}{\mathcal{V}}
\newcommand{\mcalc}{\mathrm{calc}}
\newcommand{\rD}{\mathrm{D}}
\newcommand{\rP}{\mathrm{p}}
\newcommand{\rc}{\mathrm{c}}
\newcommand{\rv}{\mathrm{v}}
\newcommand{\eye}{\mathcal{I}}
\newcommand{\Nor}{\mathcal{N}}
\paragraph{Notation.}
Symbols are expressed with the following rules.
A vector of length $d+2$ (the number of equations), with $d$ the number of dimensions, is defined by a capital letter; for instance, the vector of conservative variables of the two-dimensional system reads
$U:=[u_1,u_2,u_3,u_4]^T$ and for consistent matrix dimensions, the gradient $\nabla U := [\partial_x U, \partial_y U]^T$ is a $d(d+2)$ column vector.
Matrices and column vectors of size $d(d+2)$ are expressed with a capital calligraphic face, e.g., the mass matrix reads $\mathcal{M}$ and the $i,j$ entry is $\{\mathcal{M}\}_{ij}$.
Also, for a discrete function $w_h$ in the approximation space $V_h$, $\nrm{w_h}_{K}$ is the $L_2$ norm over $K$ and
the symbol $\pth{w}^K_h$ is a vector of length $N_p$, containing the nodal values of an element $K$; if $K$ is the reference element, the superscript is omitted.
\subsection{The compressible Navier-Stokes system}
\label{subsection:comp_ns}
The two-dimensional compressible Navier-Stokes System (cNS) in vector form consisting of the conservation of mass, momentum and energy equations, reads
\begin{equation}
	\begin{aligned}
		 & \pdv{U}{t} + \nabla \cdot \big(\mathcal{F}_\rc\pth{U} -  \vect{v}_{\mathrm{w}} \otimes U
		- \mathcal{F}_\rv\pth{U,\nabla U} \big) = 0,
	\end{aligned}
	\label{eq:cnp}
\end{equation}
with
\begin{equation*}
	\begin{aligned}
		 & U =
		\begin{bmatrix}
			\rho \\ \rho v_1 \\ \rho v_2 \\ \epsilon
		\end{bmatrix},
		 & \mathcal{F}_{\rc,x} =
		\begin{bmatrix}
			\rho v_1         \\
			\rho v_1^2 + \rP \\
			\rho v_1 v_2     \\
			v_1 \pth{\rP + \epsilon}
		\end{bmatrix},
		 &                       & \mathcal{F}_{\rc,y} =
		\begin{bmatrix}
			\rho v_2         \\
			\rho v_1 v_2     \\
			\rho v_2^2 + \rP \\
			v_2 \pth{\rP + \epsilon}
		\end{bmatrix},
		 &                       & \mathcal{F}_\rc = \big[\mathcal{F}_{\rc,x},\ \mathcal{F}_{\rc,y}\big]^T,
		 &                       &                                                                          & \mathcal{F}_\rv = \mathcal{G}(U)\nabla U,
	\end{aligned}
\end{equation*}
denoting the vector of conservative variables, the advection flux $\mathcal{F}_{\rc}$ with its $x,y$-components and the viscous flux $\mathcal{F}_{\rv}$ respectively. The divergence ($\nabla \cdot $) in this expression is applied component-wise. The $4\times 4$ block matrix $\mathcal{G}(U)$ can be expressed as a function of the primitive variables $Q = \brk{\rho, v_1, v_2, \rP}^T$ and the total energy $\epsilon = c_V \rho \vartheta + 1/2\rho \vect{v}^2$. Its entries are given by
\begin{equation*}
	\begin{aligned}
		 & \mathcal{G}_{11} = \frac{\mu}{\rho}
		\begin{pmatrix}
			0                                             & 0                          & 0 & 0 \\
			-\frac{4}{3}v_1                               & \frac{4}{3}                & 0 & 0 \\
			-v_2                                          & 0                          & 1 & 0 \\
			\bcol{\mathcal{G}_{11}}_{41}                  &
			\big(\frac{4}{3} - \frac{\gamma}{\Pr}\big)v_1 &
			\big(1 -  \frac{\gamma}{\Pr}\big)v_2          & \frac{\gamma}{\mathrm{Pr}}
		\end{pmatrix},
		\mathcal{G}_{12} = \frac{\mu}{\rho}
		\begin{pmatrix}
			0                   & 0   & 0               & 0 \\
			\frac{2}{3}v_2      & 0   & -\frac{2}{3}    & 0 \\
			-v_1                & 1   & 0               & 0 \\
			-\frac{1}{3} v_1v_2 & v_2 & -\frac{2}{3}v_1 & 0
		\end{pmatrix},
	\end{aligned}
\end{equation*}
\begin{equation*}
	\begin{aligned}
		 & \mathcal{G}_{22} = \frac{\mu}{\rho}
		\begin{pmatrix}
			0               & 0                                                    & 0                          & 0 \\
			-v_1            & 1                                                    & 0                          & 0 \\
			-\frac{4}{3}v_2 & 0                                                    & \frac{4}{3}                & 0 \\
			\bcol{\mathcal{G}_{22}}_{41}
			                & \big(1-\frac{\gamma}{\mathrm{Pr}} \big)v_1
			                & \big(\frac{4}{3}-\frac{\gamma}{\mathrm{Pr}} \big)v_2 & \frac{\gamma}{\mathrm{Pr}}
		\end{pmatrix},
		\mathcal{G}_{21} =  \frac{\mu}{\rho}
		\begin{pmatrix}
			0                  & 0               & 0   & 0 \\
			-v_2               & 0               & 1   & 0 \\
			\frac{2}{3}v_1     & - \frac{2}{3}   & 0   & 0 \\
			-\frac{1}{3}v_1v_2 & -\frac{2}{3}v_2 & v_1 & 0
		\end{pmatrix},
	\end{aligned}
\end{equation*}
whereby
\begin{equation*}
	\begin{aligned}
		 & \bcol{\mathcal{G}_{ii}}_{41} =  -\frac{1}{3}v_i^2 - \vect{v}^2 -
		\frac{\gamma}{\Pr}\Big(\frac{\epsilon}{\rho}-\vect{v}^2\Big) \qquad i=1,2.
	\end{aligned}
	\label{eq:diffusion_tensor}
\end{equation*}
The system closes with the equation of state:
\begin{align*}
	 & \rP = \pth{\gamma - 1}\rho (\epsilon - \frac{1}{2}\vect{v}^2 ),
\end{align*}
with the remaining symbols given in \autoref{tab:ns_symbols}.
\begin{table}[ht]
	\centering
	\ra{1.}
	\begin{tabular}{c|l}
		\toprule
		Symbol       & description                                \\
		\midrule
		$v_i$        & the velocity components                    \\
		$\mathrm{p}$ & pressure                                   \\
		$\rho$       & density                                    \\
		$\mu$        & dynamic viscosity                          \\
		$\gamma$     & Poisson's adiabatic constant               \\
		$\Pr$        & Prandtl number                             \\
		$\vartheta$  & Temperature                                \\
		$c_V$        & Specific heat capacity for constant volume \\
		$\mathrm{H}$ & Total enthalpy                             \\
		\bottomrule
	\end{tabular}
	\caption{Symbol description for the cNS system.}
	\label{tab:ns_symbols}
\end{table}
Regarding boundary closure, we decompose the boundary $\Gamma_b$ of the spatial computational domain $\Omega$ into solid adiabatic wall $\Gamma_{\mathrm{w}}$ and far-field (subsonic inflow and outflow) $\Gamma_{\mathrm{far}}$, so that $\Gamma_b = \Gamma_{\mathrm{w}} \cup \Gamma_{\mathrm{far}}$. On the far-field region nominal values are set, e.g., $\rP_{\infty}$, $\vect{v}_{\infty}$, and $\vartheta_{\infty}$.
On a solid adiabatic wall surface we require that $ \partial_{\vect{n}}\rP = \partial_{\vect{n}}\vartheta = 0$ and $\vect{v} = \vect{v}_{\mathrm{w}}$, with $\partial_{\vect{n}}$ denoting the derivative in the normal direction $\vect{n}$ and $\vect{v}_{\mathrm{w}}$ denoting the velocity of the solid moving body, which is needed when implementing fluid-structure interactions.

\subsection{Spatial discretisation by the discontinuous Galerkin method}
We consider an unstructured mesh $\mathcal{T}=\{K\}$, of triangles $K$. The elements are represented by a (affine or isoparametric) mapping denoted by $T_K : \hat{K} \to K$ from the reference element $\hat{K}= \bcol{(-1,-1), (1,-1), (-1,1)}$.

On $\mathcal{T}$, we consider the element-wise discontinuous vectorial polynomial approximation space $V_h$, defined as:
\begin{equation}
	\begin{gathered}
		V_h = \Big\{ \Phi_h \in \big[ L_2\pth{\Omega} \big]^4: \Phi_h\circ T_K \in \big[\mathbb{P}_p(\hat{K})\big]^4 \ \forall K \in \mathcal{T} \Big\},
	\end{gathered}
\end{equation}
with $\mathbb{P}_p(\hat{K})$ denoting the space of polynomials of degree at most $p$ over $\hat{K}$, having dimension $N_p := \dim \mathbb{P}_p(\hat{K}) = (p+1)(p+2)/2$, and $\circ$ denoting the composition of functions.

The numerical solution $U_h\in V_h$ is, thus, element-wise discontinuous and is given by the expression
\begin{equation*}
	\begin{gathered}
		U_h = \bigoplus_K U_h^K = \bigoplus_K
		\begin{cases}
			\dsum_{n=1}^{N_p} \hat{U}_{n}^K \Phi_n^K
			\pth{\vect{x}}, &
			\vect{x} \in K                       \\
			0,              & \vect{x} \not\in K
		\end{cases}, \quad  K \in \mathcal{T}.
	\end{gathered} \label{eq:Uh}
\end{equation*}

We also denote by $\Gamma:=\cup_{T\in\mathcal{T}}\partial T$ the skeleton of the mesh $\mathcal{T}$, while $\Gamma_{\mathrm{int}}:=\Gamma\backslash \Gamma_{b}$ denotes the interior part of the skeleton.
Further, using the block row vector $\mathcal{N}^{\pm} = (\vect{n}^{\pm})^T\otimes \mathcal{I}_{4\times 4}$, whereby $\vect{n} = [n_1,n_2]^T$, we define the jump
$\dbrk{\cdot}$ and the average $\dbcol{\cdot}$ operators in \autoref{tab:jump_and_average}, with $\vect{n}^\pm$ denoting the outward unit normal vector on the face $e=\partial K^+\cap \partial K^-$ shared by two elements $K^+,K^-\in\mathcal{T}$.
\begin{table}[ht]
	\centering
	\ra{1.4}
	\begin{tabular}[c]{l|l}
		\toprule
		$ e \subset \Gamma_{\mathrm{int}}$                                      & $ e \subset \Gamma_b$                     \\
		\midrule
		$\dbrk{\Phi} =  \big(\Nor^+\big)^T \Phi^+ + \big(\Nor^-\big)^T \Phi^- $ & $\dbrk{\Phi} = \big(\Nor^+\big)^T \Phi^+$ \\
		$\dbcol{\Phi} = \frac{1}{2}\pth{\Phi^+ + \Phi^-}                      $ & $\dbcol{\Phi} = \Phi^+                  $ \\
		$\dbcol{\mathcal{A}} = \frac{1}{2}\pth{\mathcal{A}^+ + \mathcal{A}^-}$  & $\dbcol{\mathcal{A}} = \mathcal{A}^+    $ \\
		\bottomrule
	\end{tabular}
	\caption{Definitions of jump and average operators.}
	\label{tab:jump_and_average}
\end{table}
The cNS system discretised by the IP DG method reads: for all $t\in I$, $I\subset \mathbb{R}$ time interval, find $U_h\equiv U_h(t)\in V_h$, such that
\begin{equation}
	\begin{aligned}
		\int_{\Omega} \pdv{U_h}{t} \cdot \Phi_h \, d\vect{x}\, & -
		\int_{\Omega} \big(\mathcal{F}_\rc\pth{U_h} -  \vect{v}_{\mathrm{w}} \otimes U - \mathcal{G}\pth{U_h}\nabla_h U_h\big) \cdot \nabla_h \Phi_h\, d\vect{x}                          \\
		                                                       & + \int_{\Gamma_{\mathrm{int}}} H \, d\vect{s} + \int_{\Gamma_b} H_b\, d\vect{s} = 0, \ \ \text{ for all } \Phi_h\in V_h,
	\end{aligned}
	\label{eq:weak_form}
\end{equation}
with
\begin{align*}
	 & H = H_\rc - H_\rv, \ \  H_{\rc}\equiv  H_{\rc}\pth{U_h^\pm,\Phi_h^\pm;\vect{n}},\ \ H_{\rv}\equiv H_{\rv}\pth{U_h^\pm,\Phi_h^\pm;\vect{n}}
\end{align*}
the interior face terms containing the numerical fluxes for advection (subscript $\rc$) and diffusion (subscript $\rv$) and $H_b$ comprising all boundary-related terms; these will be given precisely below.

\subsubsection{Advection numerical flux}
The advection flux $H_\rc$ includes a ``central'' averaging term and an \say{upwind} term enforcing numerical diffusion in accordance to the system's local advection
eigenvalues. In particular, we use the Roe flux \cite{roe1981approximate} as incorporated in the DG setting in \cite{persson2009discontinuous}. Roe's numerical flux is constructed so that the, so-called, \textit{Rankine-Hugoniot} property
\begin{equation*}
	\mathcal{F}^+_{\mathcal{N}} - \mathcal{F}^-_{\mathcal{N}} = \mathcal{A}_{\text{visc}} (U_h^+ - U_h^-),
\end{equation*}
whereby $\mathcal{F}^{\pm}_{\Nor} = \Nor^{\pm} \big(\mathcal{F}_c \pth{U^{\pm}_h} - \vect{v}_{\mathrm{w}} \otimes  U^{\pm}_h\big)$,
holds at the discrete level. The numerical flux reads
\begin{equation}
	H_\rc = \big(  \dbcol{\mathcal{F}_{\Nor}}  - \dbrk{\mathcal{A}_{\text{visc}} U_h} \big)\cdot \dbrk{\Phi_h},
	\label{eq:roe_flux}
\end{equation}
with $\mathcal{A}_{\text{visc}}$ the Roe matrix diagonalised as $\mathcal{A}_{\text{visc}} = \mathcal{R}^{-1} \Lambda \mathcal{R}$ and the eigenvalues forming the diagonal $\Lambda$ being $\{v_{\alpha},\ v_{\alpha} \pm c\}$, whereby
\begin{align*}
	 & v_{\alpha} = |(\Tilde{\vect{v}} - \vect{v}_{\mathrm{w}})\cdot \vect{n}^{\pm}|, &                                          & c^2 = (\gamma-1)\big(\Tilde{\mathrm{H}} - \frac{1}{2}\Tilde{\vect{v}}^2\big),
	 &                                                                                & \mathrm{H} = \epsilon + \mathrm{p}/\rho,
\end{align*}
and $(\Tilde{\cdot})$ used here to express a Roe-averaged value.
Other choices include the Vijayasundaram flux \cite{hartmann2008optimal}, the Lax-Friedrichs flux \cite{winters2018comparative}, and the HLLC \cite{cheng2016direct}.
\subsubsection{Diffusion numerical flux}
The diffusion flux stems from the interior-penalty flux from \cite{hartmann2008optimal}, reading
\begin{equation}
	\begin{aligned}
		 & H_\rv = \dbcol{\mathcal{G}(U_h)\nabla U_h} \cdot \dbrk{\Phi_h}\, +\,
		\theta \dbcol{\mathcal{G}(U_h)\nabla \Phi_h} \cdot \dbrk{U_h}\, -\, \sigma_{e} \dbrk{U_h} \cdot \dbrk{\Phi_h},
	\end{aligned}
	\label{eq:ip_flux}
\end{equation}
with the first term on the right-hand side ensuring consistency,
the second term symmetrising the first term in a consistent fashion, while the last term, is typically referred to as the \emph{discontinuity-penalisation} or, simply, \emph{penalty} term with parameter $\sigma_e\ge 0$; its function is to penalise the jump across element boundaries.
The parameter $\theta\in[-1,1]$ defines different variants: the non-symmetric ($\theta=-1$), the incomplete ($\theta=0$), or the symmetric $(\theta=1$) interior penalty method, respectively.

The penalty parameter $\sigma_e\ge 0$ has to be chosen sufficiently large face-wise to ensure positivity of the discretisation of the diffusion. At the same time, as we shall see below, large values of $\sigma_e$ impose a more stringent CFL condition. Therefore, accurate estimation of $\sigma_e$ is extremely important for the practical implementation of the method. To that end, a new, refined choice of $\sigma_e$, sensitive to flow characteristics and the nonlinearity of the diffusion, is presented below, allowing for less stringent CFL constants than other recipes for $\sigma_e$ from the literature.

\subsubsection{The interior penalty term}

The penalty term is required for the coercivity/positivity of the diffusion part of the operator. To see this, we begin by recalling the following trace inverse estimate \cite{warburton2003constants}: for $v\in \mathbb{P}_p(K)$, and for each face $e\subset \partial K$, we have
\begin{equation}
	\nrm{v}_{e}^2 \leq C_{\rm inv}(e,K,p)\nrm{v}_{K}^2,\quad   \text{ for }\quad    C_{\rm inv}(e,K,p):=\frac{\pth{p+1}\pth{p+d}|e|}{d|K|}.
	\label{eq:trace_inverse_estimate}
\end{equation}
Then, writing $\mathcal{G}\equiv\mathcal{G}(U_h)$, we set
$
	\bar{\mathcal{G}}_e:=\max_{*\in\{+,-\}}
	\|\,
	| \mathcal{N}\mathcal{G}|_{K^*}\mathcal{N}^T |_2
	\,\|_{L_\infty(e)}
$,
for brevity, with $e=\partial K^+\cap\partial K^-$ a generic interior face. Using, now H\"older/Cauchy-Scwharz inequalities along with
\eqref{eq:trace_inverse_estimate}, we have, respectively,
\begin{equation*}
	\begin{aligned}
		\int_{\Gamma_{\text{int}}}\dbcol{\mathcal{G}\nabla U_h}\cdot\dbrk{\Phi_h} \, d\vect{s}\,
		\le\, & \ \frac{1}{2}\sum_{e\subset \Gamma_{\text{int}}}
		\bar{\mathcal{G}}_e\Big(\sum_{*\in\{+,-\}}\!\nrm{\tau_e^{-\frac{1}{2}} \nabla U_h^*}_e\Big) \nrm{\sqrt{\tau_e} \dbrk{\Phi_h}}_e \\
		\le\, & \ \frac{1}{2}\sum_{e\subset \Gamma_{\text{int}}}
		\Big(\frac{\bar{\mathcal{G}}_e^2}{\tau_e}\max_{*\in\{+,-\}}\!C_{\rm inv}(e,K^*,p)\Big)^{\frac{1}{2}}\Big(\sum_{*\in\{+,-\}} \!\nrm{\nabla U_h}_{K^*}\Big)\nrm{\sqrt{\tau_e}\dbrk{\Phi_h}}_e.
	\end{aligned}
\end{equation*}
Selecting, now,
\begin{equation*}
	\tau_e:=(1+\theta)^2(d+1)\bar{\mathcal{G}}_e\max_{*\in\{+,-\}}C_{\rm inv}(e,K^*,p),\quad\text{ for }\ \theta\in(-1,1],
\end{equation*}
we deduce
\begin{equation}
	\begin{aligned}
		(1+\theta)\int_{\Gamma_{\text{int}}}\dbcol{\mathcal{G}\nabla U_h}\cdot\dbrk{\Phi_h} \, d\vect{s}
		\ \le & \ \frac{1}{2}\sum_{e\subset \Gamma_{\text{int}}}
		\Big(\sum_{*\in\{+,-\}} \!\Big(\frac{\bar{\mathcal{G}}_e}{d+1}\Big)^{\frac{1}{2}} \nrm{\nabla U_h}_{K^*}\Big) \nrm{\sqrt{\tau_e}\dbrk{\Phi_h}}_e                                                       \\
		\ \le & \ \frac{1}{2}\sum_{e\subset \Gamma_{\text{int}}} \sum_{*\in\{+,-\}} \frac{\bar{\mathcal{G}}_e}{d+1} \nrm{\nabla U_h}_{K^*}^2+\frac{1}{2} \nrm{\sqrt{\tau_e}\dbrk{\Phi_h}}_{\Gamma_{\rm int}}^2 \\
		\ \le & \ \frac{1}{2}\nrm{\sqrt{\bar{\mathcal{G}}|_{K}}\nabla_h U_h}_{\Omega}^2+ \frac{1}{2} \nrm{\sqrt{\tau_e}\dbrk{\Phi_h}}_{\Gamma_{\rm int}}^2,
	\end{aligned}
	\label{inverse_magic}
\end{equation}
where $\bar{\mathcal{G}}|_K:=\max_{e\subset\partial K}\bar{\mathcal{G}}_e$.
Hence, the penalty term reads:
\begin{equation}
	\sigma_e \equiv C_1 \tau_e= C_1 (1+\theta)^2(d+1)\bar{\mathcal{G}}_e \max_{*\in\{+,-\}}C_{\text{inv}}(e,K^*,p),
	\label{eq:penalty_final}
\end{equation}
where $C_1\ge 0$ is a free parameter that has to be chosen large enough for stability, at least for implicit methods. Its values, however, are expected to be insensitive to $\theta,p,K^*,e,\bar{\mathcal{G}}_e$ due to the above estimation.

For comparison, standard choices of penalty parameter in the literature are of the form $\sigma|_e=C_\sigma p^2/\min\{h_{K^-},h_{K^+}\}$, on each face $e=\partial K^-\cap \partial K^+$, for $K^-,K^+\in\mathcal{T}$ neighbouring elements \cite{hartmann2006adaptive, hartmann2008optimal, zou2025moving}, with the constant $C_\sigma$ typically selected as $C_{\sigma}\ge 10$ or larger in case of irregular/anisotropic elements. In contrast, the constant $C_1$ in \eqref{eq:penalty_final} is oblivious to mesh anisotropy. Overall, the choice \eqref{eq:penalty_final} yields typically smaller penalty constants than the standard choice $C_{\sigma}$ from the literature.

\subsubsection{The boundary term}
In order to enforce boundary conditions (BC) on $e\subset \Gamma_b$, a \textit{ghost} cell $K^b \not \subset \Omega$ is considered, for which $e=\partial K^+\cap \partial K^b,\ K^+ \subset \Omega$.
Here, two types of BC are needed; these are, the inflow-outflow and the \say{moving} solid adiabatic wall conditions.

On $\Gamma_{\mathrm{far}}$, the boundary values are $U^b = \brk{\rho_{\infty}, \rho_{\infty} \vect{v}_{\infty}, \epsilon_{\infty}}^T$, the normal advection boundary flux is approximated by $H_\rc(U^+,U^b, \Phi^+;\vect{n}^+)$  and the normal diffusion boundary flux is $H_{\rv} (U_h^+,U^b,\Phi^+;\vect{n}^+)$. The total far-field boundary flux reads:
\begin{align*}
	 & \left. H_b\right|_{e \subset \Gamma_{\mathrm{far}}} = H_c\pth{U_h^+,U^b}
	- \dbcol{ \mathcal{G}\pth{U_h}\nabla U_h} \cdot \dbrk{\Phi_h}
	- \theta \dbcol{ \mathcal{G}\pth{U_h}\nabla \Phi_h} \cdot \dbrk{U_h}
	+ \sigma_e\, (\Nor^+)^T (U_h^+ - U^b)\cdot \dbrk{\Phi_h}
	.
\end{align*}
On $\Gamma_{\mathrm{w}}$, Neumann-type BC are needed for $\rP$ and $\vartheta$ and Dirichlet-type BC for $\vect{v}$. Thus, the boundary primitive vector is $Q^b = \brk{\rho^+, \vect{v}_{\mathrm{w}}, \mathrm{p}^+}^T$, from which $U^b$ is calculated. In this occasion, the total wall boundary flux reads:
\begin{align*}
	 & \left. H_b\right|_{e \subset \Gamma_{\mathrm{w}}} = \Nor^+ \big(\mathcal{F}_c\pth{U^b} - \vect{v}_{\mathrm{w}} \otimes U^{b} \big)
	- \mathcal{G}\pth{U^b} \nabla U_h^+ \cdot \dbrk{\Phi_h}
	- \theta\, \mathcal{G}\pth{U^b}\nabla \Phi^+_h \cdot \dbrk{U_h}
	+ \sigma_e\, (\Nor^+)^T (U_h^+ - U^b) \cdot \dbrk{\Phi_h}.
\end{align*}
In the following, the nodal formulation is briefly presented.
\subsubsection{Weak matrices and nodal formulation}
\begin{figure}
	\centering
	\includegraphics[width=0.3\linewidth]{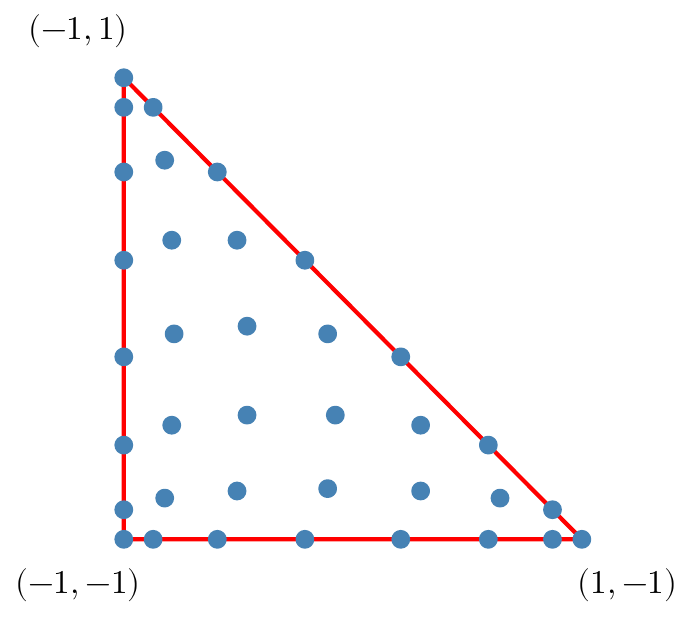}
	\caption{Distribution of interpolating nodes on $\hat{K}$ for $p=5$.}
	\label{fig:optimal_points}
\end{figure}

In the nodal approach, the degrees of freedom (dof) are defined as interpolating values on $\hat{K}$ positioned at carefully-chosen interpolation cites, which are evaluated by utilising the procedure presented in \cite{hesthaven2008nodal}. Their distribution for $p=5$ is shown in \autoref{fig:optimal_points}.

If we set $\mathcal{F}:=\mathcal{F}_\rc - \mathcal{F_\rv}$, $H:=H_{\rc} - H_{\rv}$ and choose the Lagrange base
$\ell_i^K\equiv \ell_i^K(\vect{x}),\ i=1,\dotsc,N_p$, $K\in \{\mathcal{T}\}$, then the semi-discrete expression for element $K$ and for $k=1,\dotsc,4$ reads:
\begin{equation*}
	\int_K \pdv{\{U_h^K\}_k}{t} \ell _h^K\, d\vect{x}\, - \int_K \{\mathcal{F}\pth{U_h^K}\}_k \cdot \nabla \ell_h^K\, d\vect{x}\, +
	\int_{\partial K} \{H\}_k\, d\vect{s} = 0.
\end{equation*}

If, for $K\equiv K^+$ with $e= K^+ \cap K^-$, we set $\Nor^- = -\Nor^+$, then:
\begin{align*}
	\{H^K_{\mathrm{m}}\}_k & = \big\{ \dbcol{\mathcal{F}_\rc - \mathcal{G}(U_h)\nabla U_h}  - (\mathcal{A}_{\text{visc}} -\sigma_e)\Nor^+ \dbrk{U_h} \big\}_k\ \ell_i^e, &  &                    \\
	\{H^K_{\theta,j}\}_k   & = - \theta\big\{(U_h^+ - U_h^-)^T \Nor^+ \mathcal{G}\pth{U^+_h}\big\}_{4j+k-4}\, \partial_{x_j} \ell_i^K,                                                          &  & \text{for}\ j=1,2,
\end{align*}
and, of course, $H  = H_{\mathrm{m}} + H_{\theta,1} + H_{\theta,2}$.

The mass and stiffness \textit{element-local} matrices read:
\begin{equation}
	\begin{aligned}
		\big\{ \mathcal{M}^K \big\}_{ij}
		 & = \int_K \ell^K_i\ell^K_j\, d\vect{x},
		 &                                                                 & \big\{ \pth{\mathcal{M} \{H_{\mathrm{m}}\}_k}_h^e\big\}_i  = \int_e \{H_{\mathrm{m}}\}_k \, \ell^e_i\, d\vect{x},            \\
		\big\{\pth{\mathcal{S}_x \{\mathcal{F}\}_k}^K_h\big\}_i
		 & = \int_K \{\mathcal{F}\}_k\, \partial_x \ell^K_i\, d\vect{x},
		 &                                                                 & \big\{ \big( \mathcal{S}_y \{\mathcal{F}\}_{k+4}\big)^K_h\big\}_i =
		\int_K \{\mathcal{F}\}_{k+4} \partial_y\ell_i^K\, d\vect{x},                                                                                                                                      \\
		\big\{\pth{\mathcal{S}_x \{H_{\theta x}\}_k}_h^K\big\}_i
		 & = \int_e \{H_{\theta x}\}_k \, \partial_x \ell_i^K\, d\vect{x},
		 &                                                                 & \big\{\pth{\mathcal{S}_y \{H_{\theta y}\}_k }_h^K\big\}_i    = \int_e \{H_{\theta y}\}_k \, \partial_y \ell_i^K\, d\vect{x}.
		\label{eq:weak_matrices_calc}
	\end{aligned}
\end{equation}

We note that, unlike the mass matrix $\mathcal{M}^K$ which is defined as a dot product of the polynomial base, the stiffness matrices $\mathcal{S}_x,\mathcal{S}_y$ contain the flux functions which have non-linear components with respect to $U$. In fact, many entries of $\mathcal{F}(U_h)$ and $H(U_h)$ are rational functions and the underestimation of these terms can lead to aliasing-driven instabilities. These are thoroughly studied for DG methods in \cite{kirby2003aliasing,Gassner2013221,spiegel2015aliasing,winters2018comparative} where the cNS is considered with subsonic initial conditions. In \cite{Gassner2013221} it is proposed to address this by oversampling the flux functions when performing the spatial integration. This is achieved either by integrating exactly nodal interpolants on the nodes of order $p_f=3p$ basis, or by employing $\mathbb{P}_{p_f}(K)$-exact  quadrature. In the latter, often preferred, case, the symmetric cubature rules of \cite{witherden2015identification} are used.

The calculation of the local matrices is performed on the reference element. Since the Jacobian of $T_K$ is constant on $K$, these are evaluated and stored only for $\hat{K}$ that offers significant savings in memory usage. Further details can be found in \cite{hesthaven2008nodal}.

Using the expressions from \eqref{eq:weak_matrices_calc}, the semi-discrete expression from \eqref{eq:weak_form} for every $K\in \{\mathcal{T}\}$ and all components $k=1,\dotsc,4$ reads $\mathcal{M}^K \pdv{\pth{\{U\}_k}^K_h}{t} = \pth{\{R\}_k}_h^K$, with:
\begin{equation}
	\begin{aligned}
		\pth{\{R\}_k}_h^K = &
		\big(\mathcal{S}_x \{\mathcal{F}\}_{k}\big)_h^K
		+ \big(\mathcal{S}_y \{\mathcal{F}\}_{k+4}\big)_h^K
		- \big(\mathcal{S}_x \{H_{\theta x}\}_k \big)_h^K
		- \big(\mathcal{S}_y \{H_{\theta y}\}_k \big)_h^K
		- \sum_{e \subset \partial K} \mathcal{E}^{e} \pth{\mathcal{M}   \{H_{\text{m}}\}_k }_h^e
	\end{aligned}
	\label{eq:dg_on_K}
\end{equation}
being the residual and $\mathcal{E}^e$ representing a connectivity mapping of size $N_p \times (p+1)$ that maps the nodal points attached to face $e$ to the element-local numbering of $K$.

\subsection{Time discretisation and the CFL plateau}

\label{subsection:time_discretisation}

In the context of CFD generally, explicit time discretization is only used in cases where the
fully-discrete system is \textit{globally} non-stiff. The discrete system's stiffness becomes profound in cases of external boundary
layer flows where the element size spans several length scales or when inflow/outflow boundary conditions are close to the incompressible limit.

In this study, an explicit time marching was found to be sufficient for laminar 2D VIV test cases; hence, a low storage Strong-Stability-Preserving (SSP) Runge-Kutta scheme of $N_{\text{RK}}=5$ stages from \cite{niegemann2012efficient} is used.
The implementation of an RKDG method carries the assumption that the solution is \say{element-explicit}, thus requiring the inversion of the mass matrix. However, since the global mass matrix is block diagonal, the inversion is direct.

For a grid $\{t_n\} \subset I,\ n\in [0,N]$, let $\vect{U}_n\equiv U_h(t_n)$ be the unknown vector, $\vect{R}_n$ the residual, and $\vect{\mathcal{M}}$ the block-sparse mass matrix. Then, a low storage RK can be implemented as shown in Algorithm \ref{alg:low_storage_rk}, where $A_i,B_i$,$0\le i\le N_{\text{RK}}$, are the coefficients for an $N_{\text{RK}}$-stage variant (see \cite{niegemann2012efficient}).
\begin{algorithm}
    \caption{Low Storage version of a RK step}
    \setstretch{1.0}
    \label{alg:low_storage_rk}
    \begin{algorithmic}
        \State $\vect{K}_0     \gets \vect{U}_n$
        \State $\vect{K}_1     \gets A_i \vect{K}_1 + \delta t \vect{\mathcal{M}}^{-1} \vect{R}\pth{t_i + c_i \delta t, \vect{K}_0}$
        \State $\vect{K}_0     \gets \vect{K}_0 + B_i \vect{K}_1$
	\State $\vect{U}_{n+1} \gets \vect{K}_0$
    \end{algorithmic}
\end{algorithm}

Explicit time integration comes with the burden of time-step restriction imposed by the CFL condition. In Appendix \ref{CFL_Heuristics}, we present some bounds aiming to justify the following choice for the CFL condition
\begin{equation}
    \begin{aligned}
        & C_{\mathrm{CFL}} = \Lambda(\mathcal{G}, \vect{b}_{\text{ns}}, \mathcal{T}, p)\delta t, & \{\vect{b}_{\text{ns}}\}_j = |\{\vect{v} - \vect{v}_{\mathrm{w}} \}_j | + 
        \sqrt{\frac{\gamma \mathrm{p}}{\rho}},
        \label{eq:cfl_general}
        \end{aligned}
    \end{equation}
where $\delta t$ is the time-step and $\Lambda(\mathcal{G}, \vect{b}_{\text{ns}}, \mathcal{T}, p)$ as in \eqref{rayleigh_bound} in the Appendix \ref{CFL_Heuristics}.

\subsection{Assembly and coupling with rigid dynamics}

\begin{figure}
    \centering
    \includegraphics[width=0.6\linewidth]{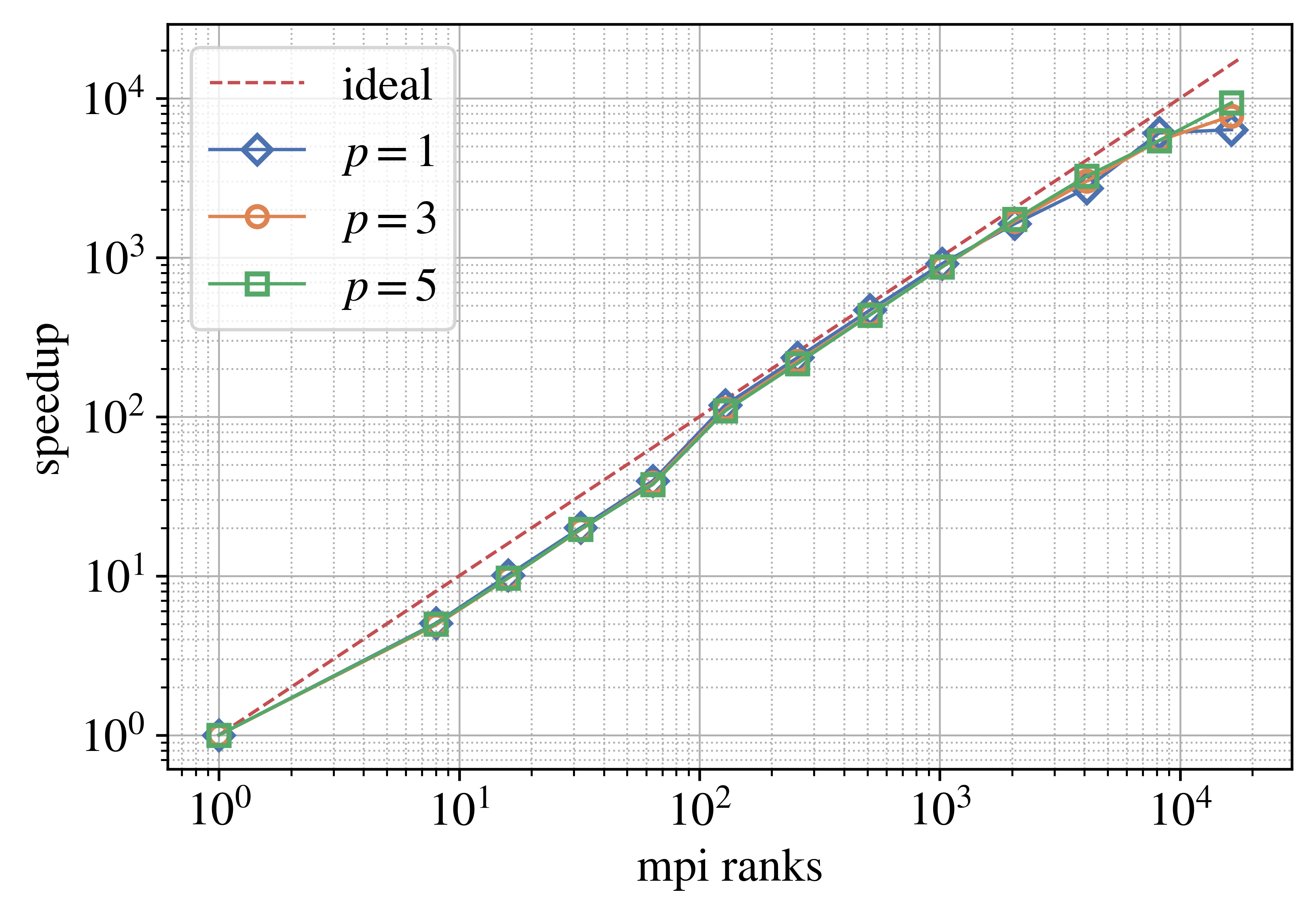}
    \caption{Strong scaling of the nodal IP RKDG algorithm. Mesh of $500,000$ cells.}
    \label{fig:strong_scaling}
\end{figure}

The element-block residual $\vect{R}_n$ is calculated in a parallel framework as shown in Algorithm \ref{alg:cfd_res}.
{\color{red}
\begin{algorithm}
    \caption{Calculation of the right hand side for the nodal IP RKDG}
    \setstretch{1.2}
    \label{alg:cfd_res}
    \begin{algorithmic}[1]
        \State  \Call{send mpi neighbour data}{}$\bcol{\pth{U}_h^e,\ \pth{\nabla U}^e_h  }$
        \State $\pth{\mathcal{S}_x\{\mathcal{F}\}_k}^K_h + \pth{\mathcal{S}_y\{\mathcal{F}\}_{k+4}}^K_h  \gets$  \Call{calculate volume integrals}{}$\pth{U^K_h}$
        \State $(\{R\}_k)_h^K \gets \pth{\mathcal{S}_x\{\mathcal{F}\}_k}^K_h + \pth{\mathcal{S}_y\{\mathcal{F}\}_{k+4}}^K_h$ 
        \State $ \pth{U}^e_h,\ e\in \Gamma_{b} \gets$  \Call{set boundary conditions}{}$(\vect{v}_{\infty},\mathrm{p}_{\infty},\theta_{\infty},\vect{v}_{\mathrm{w}})$
        \State $\pth{U}_h^e,\ \pth{\nabla U}^e_h \gets$ \Call{receive mpi neighbour data}{}
        \State \Call{calculate numerical fluxes}{}$\pth{U_h^+,U_h^-,e}$
        \State $(\{R\}_k)_h^K \gets (\{R\}_k)_h^K - \sum_{e\in K} \mathcal{E}^e \pth{\mathcal{M}_{1\rD}\{H\}_k  }_h^e
         - \pth{\mathcal{S}_x \{H_{\theta,x}\}_k}_h^K - \pth{\mathcal{S}_y \{H_{\theta,y}\}_k }_h^K$
    \end{algorithmic}
\end{algorithm}}
This algorithmic implementation reveals that parallel communication can be overlapped with \textit{cpu-local} floating point operations using asynchronous messages, which, together with the stencil compactness offered by the nodal IP RKDG can leverage the code scalability at over $15,000$ processes as shown in \autoref{fig:strong_scaling}.

The last part involves the coupling of IP RKDG with rigid Newtonian dynamics. For the current study, lateral oscillation is considered, and thus the oscillator's movement is dictated by the scalar ordinary differential equation:
\begin{equation}
    M_{\text{r}}\, \ddot{y} + C_{\text{r}}\,\dot{y} + K_{\text{r}}\, y = F_{\text{r}},
    \label{eq:1dof_oscillator}
\end{equation}
which is numerically integrated using a Newmark-$\beta$ method with $\beta_N=0.25$ and $\gamma_N=0.5$ Newmark coefficients.

Equation \eqref{eq:1dof_oscillator} is decoupled with every IP RKDG iteration, as shown in the algorithm \ref{alg:cfd_rbd}.
In detail, the grid movement velocity $\vect{v}_{\text{w}}\equiv \vect{v}_{\text{w}}(t_n) = 
\delta \vect{x}_b / \delta t$ 
\footnote{$\delta \vect{x}_b \equiv \vect{x}_b(t_{n+1}) - \vect{x}_b(t_n)$ is the translation of the solid's center of gravity.}
is retrieved as a result of a single Newmark iteration for $\delta t$ time-step dictated by the CFL condition \eqref{eq:cfl_general}. That being said, $\vect{v}_{\text{w}}$ is fixed during the RK iteration, past which the local force ($\vect{F}_{\text{local}}$) and moment ($\vect{M}_{\text{local}}$) vectors, evaluated on the partitioned parts of $\Gamma_{\text{w}}$, are summed and broadcast across all computing processes.
\begin{algorithm}
   \caption{IP RKDG with rigid body dynamics}
   \label{alg:cfd_rbd}
   \begin{algorithmic}[1]
		\State $F_{\mathrm{ext}} \gets$  \Call{Calculate Aerodynamic Forces}{}$\pth{\vect{U}_{\mathrm{init}}}$
	   \State $t \gets 0$
	   \While{$t < t_{final}$}
		\State $\delta t \gets  $         \Call{Calculate Time Step}{}$\pth{\mathrm{CFL}}$
		\State $\delta \vect{x}_b \gets $   \Call{Rigid Body Dynamics Newmark}{}$\pth{F_{\mathrm{global}}}$
		\State $\vect{v}_b \gets \delta \vect{x}_b / \delta t$
		\State $\mathcal{T} \gets$			 \Call{Update Mesh Coordinates}{}$\pth{\mathcal{T}, \vect{v}_b}$
		\State $\vect{U} \gets$	          \Call{Low Storage SSP-RK}{}$\pth{\vect{U}, \vect{v}_b,\delta t}$
		\State $(\vect{F}_{\text{local}}, \vect{M}_{\text{local}}) \gets$   \Call{Calculate Aerodynamic Forces}{}$\pth{\vect{U}}$
		\State $(\vect{F}_{\text{global}},\vect{M}_{\text{global}}) \gets$ \Call{Reduce and Broadcast Forces}{}$\pth{\vect{F}_{\text{local}}, \vect{M}_{\text{local}}}$
		\State $t \gets t + \delta t$
	   \EndWhile
   \end{algorithmic}
\end{algorithm}

\section{Validation and results}
In this section, we present the results of this study. First, a convergence test is carried out showing the accuracy order of the nodal IP RKDG framework. Then, for the von-K\'{a}rm\'{a}n vortex street behind a stationary circular cylinder, we conduct a sensitivity study and we validate with experimental and numerical data from the literature. In the last part, the application of the nodal IP RKDG framework on a VIV problem is studied and compared against available numerical results from the literature. For the rest of the study,  the \say{$s$} after an interval $I=(a,b)$ refers to \textit{seconds}.
\subsection{Convergence study}

\label{subsection:convergence_test}

This case is inspired by corresponding convergence studies in \cite{gassner2008discontinuous,hartmann2008optimal}. It entails comparing the IP RKDG solution to a manufactured known exact solution, which results to the presence of a source term.
The parameters for the exact solution are $\mu = 0.01,\ \kappa =25,\ C_2=200$,
whereby $C_2$ is chosen to correspond to a low Mach number, i.e., $\mathrm{Ma} \approx 0.15$.
The exact solution reads:
\begin{align}
    & \phi = \kappa \pth{x+y},
    & U = \big[\sin\phi + C_2,\, \sin\phi + C_2,\, \sin\phi + C_2,\, \pth{ \sin\phi + C_2}^2\big]^T,
    \label{eq:analytical_solution}
\end{align}
corresponding to the source term:
\begin{align*}
    & \gamma_1 = \gamma -1,
    & S = 
    \begin{bmatrix}
        2\kappa\cos\phi \\
        \kappa\gamma_1\sin2\phi +\pth{2\kappa - \frac{1}{2}k\gamma_1 + 2\kappa C_2\gamma_1} \cos\phi \\
        \kappa\gamma_1\sin2\phi +\pth{2\kappa - \frac{1}{2}k\gamma_1 + 2\kappa C_2\gamma_1} \cos\phi \\
        2\kappa\gamma \sin2\phi + \pth{4k\gamma C_2 - \kappa\gamma_1} \cos\phi + \frac{2\gamma \mu \kappa^2}{\Pr}  \sin\phi
    \end{bmatrix}.
\end{align*}

For the source term, overintegration is applied to decrease aliasing errors below the measurable tolerance. To this end, we use $p_f=25$.

\begin{figure}
    \centering
    \begin{subfigure}[b]{0.45\textwidth}
        \centering
        \includegraphics[width=0.8\linewidth]{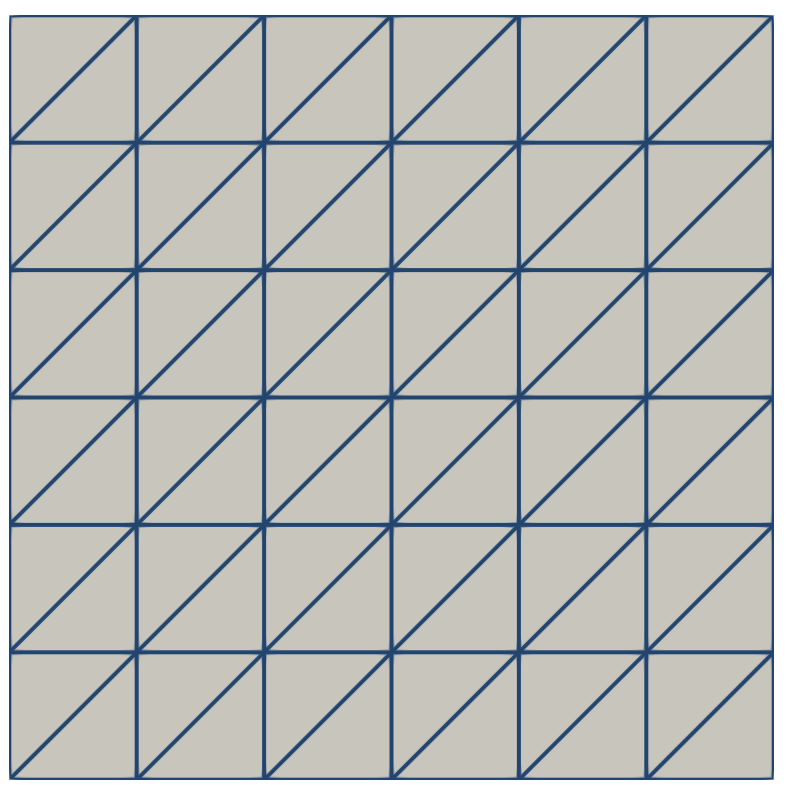}
    \end{subfigure}
    \begin{subfigure}[b]{0.45\textwidth}
        \centering
        \includegraphics[width=0.8\linewidth]{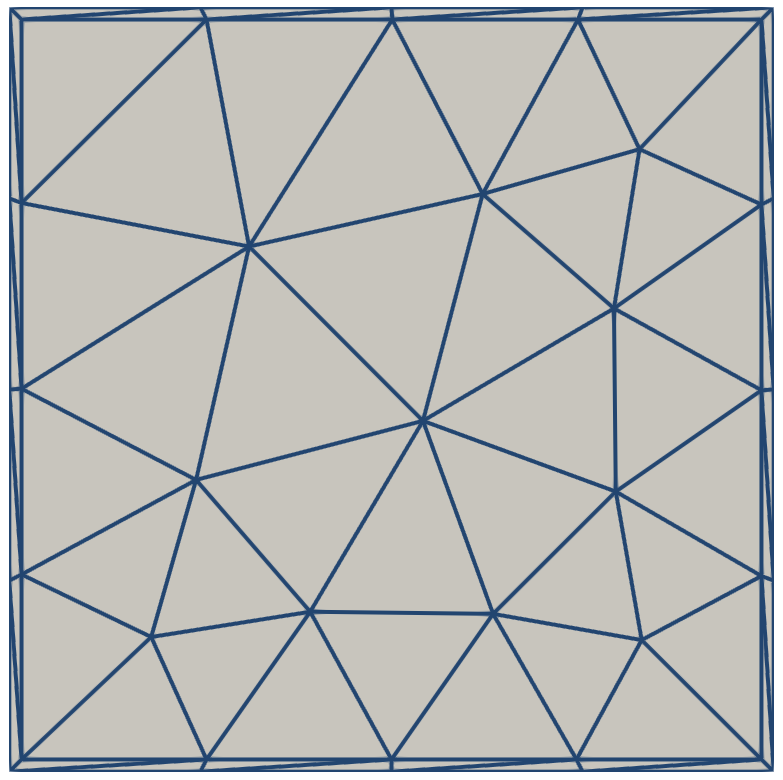}
    \end{subfigure}
    \caption{Convergence test: $h_1$ version of $m_1$ (left) and $m_2$ (right) meshes.}
    \label{fig:m1m2}
\end{figure}

The convergence test consists of approximating $U$ using both the symmetric and the incomplete IP  methods (SIP and IIP for short). For these, the penalty constant of $C_1=0.01$ is applied for both variants of interior-penalty DG. The test is carried out for $1\le p\le 5$ and three mesh refinements; starting from $h_1$ to $h_2, h_3$ and $h_4$. For all cases, we evaluate the convergence rate $r\equiv r(p)$ using linear regression. The procedure is repeated for two types of meshes: $m_1$ which has ordered cells perpendicular to the (exact) solution's direction and $m_2$ which is unstructured with one layer of flattened cells; as seen in \autoref{fig:m1m2}. 
\begin{figure}
    \centering
    \begin{subfigure}[b]{0.49\textwidth}
        \centering
        \includegraphics[width=0.99\linewidth]{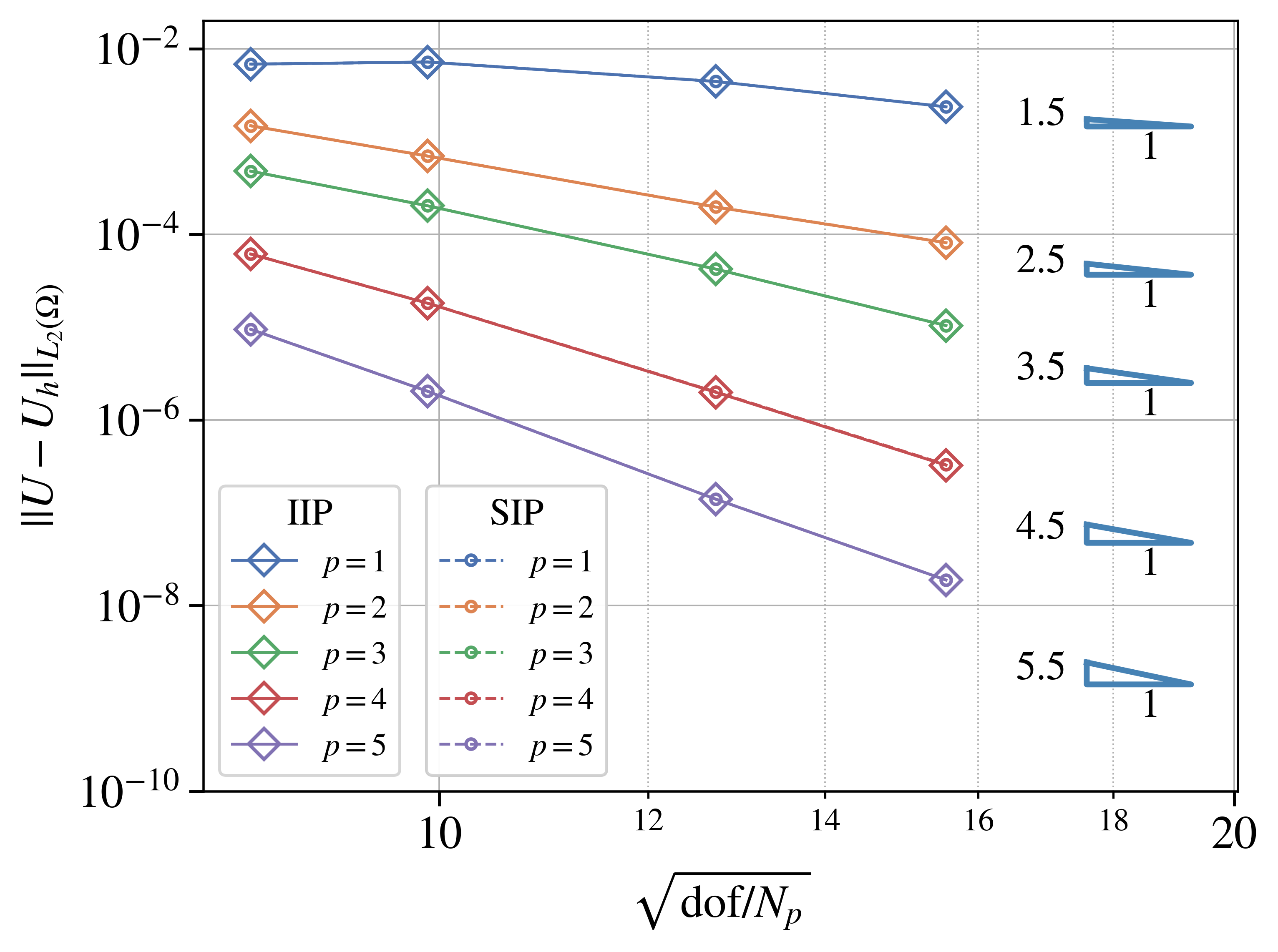}
    \end{subfigure}
    \begin{subfigure}[b]{0.49\textwidth}
        \centering
        \includegraphics[width=0.99\linewidth]{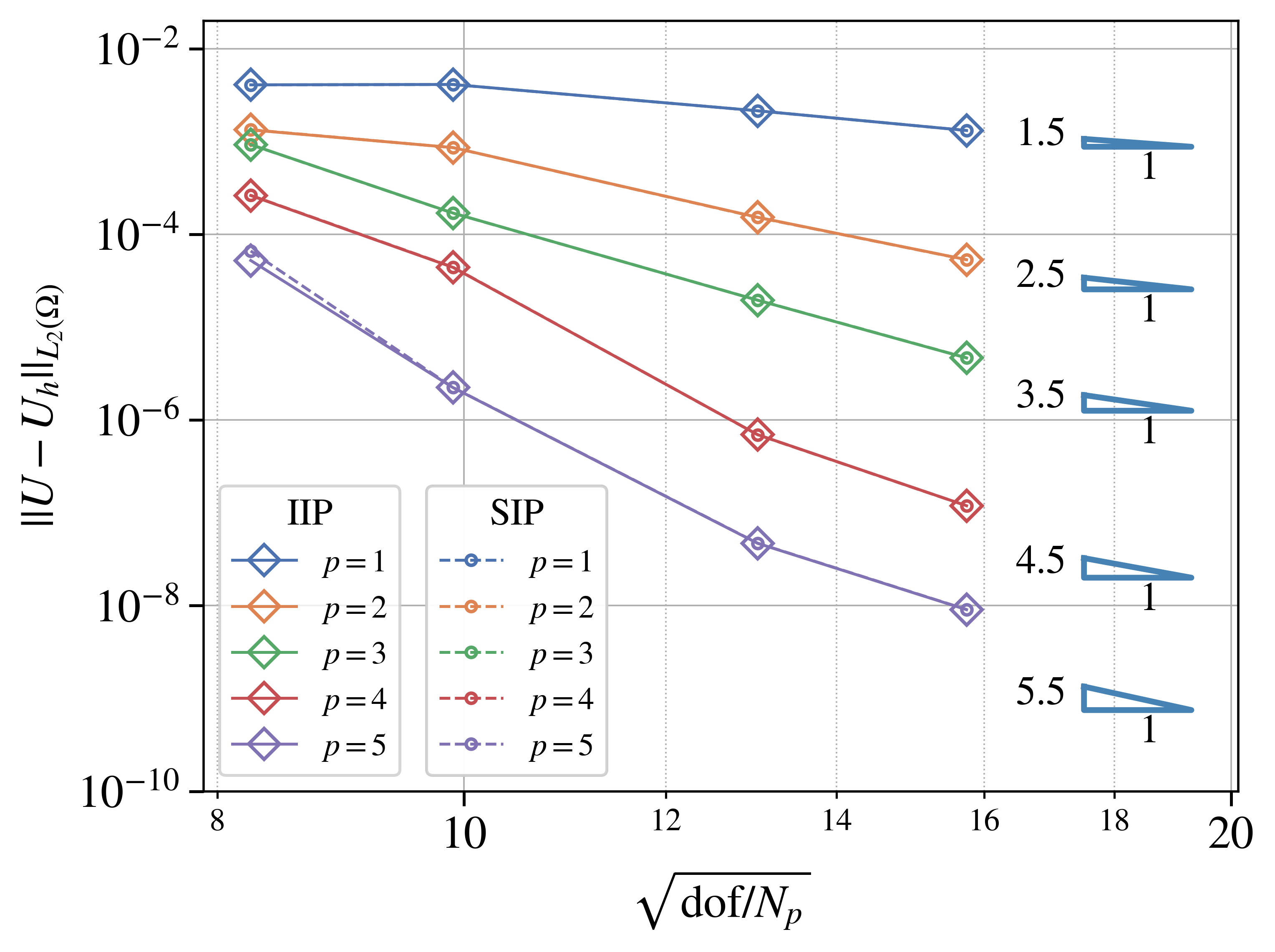}
    \end{subfigure}
    \begin{subfigure}[b]{0.49\textwidth}
        \centering
        \includegraphics[width=0.99\linewidth]{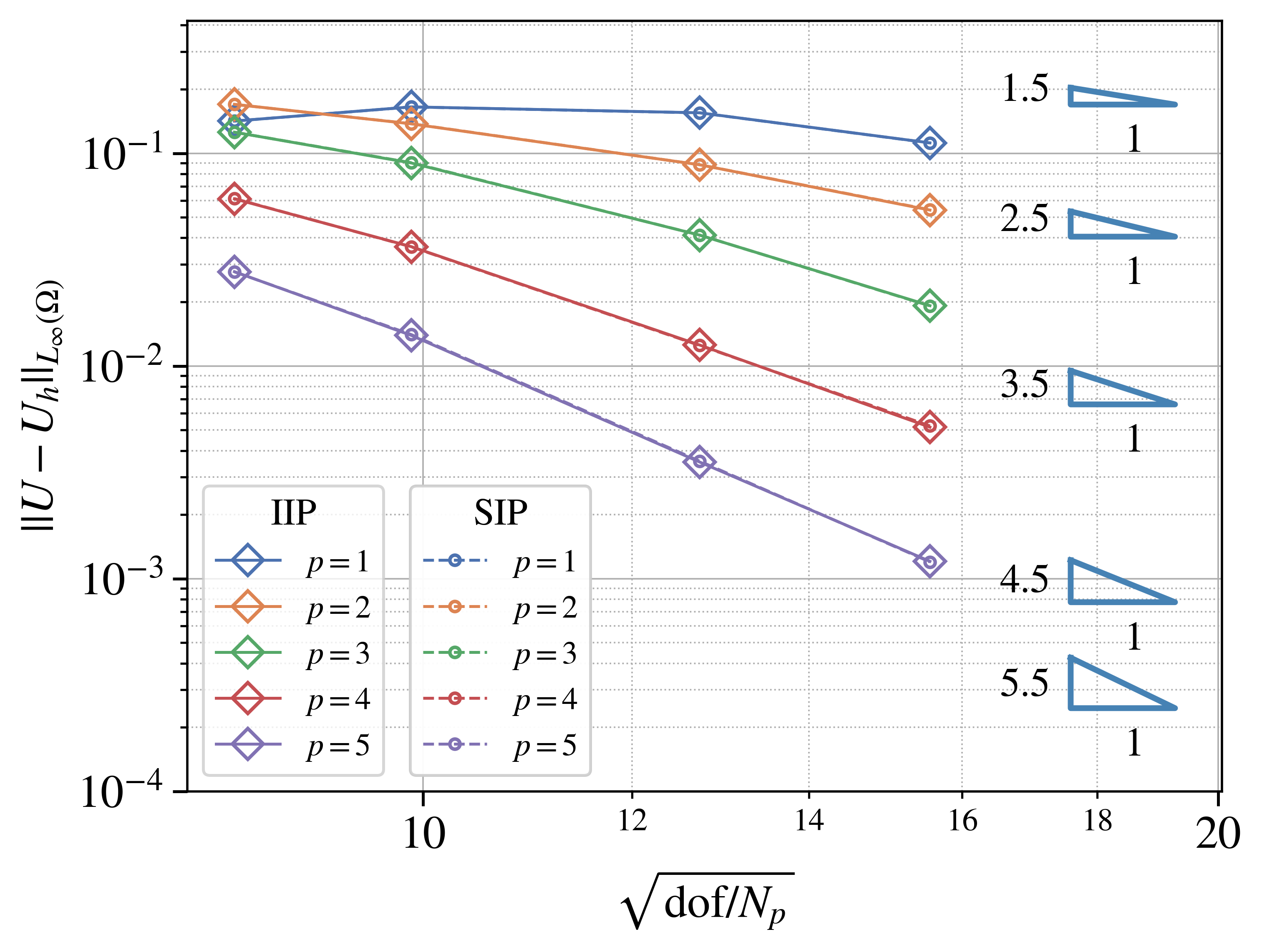}
    \end{subfigure}
    \begin{subfigure}[b]{0.49\textwidth}
        \centering
        \includegraphics[width=0.99\linewidth]{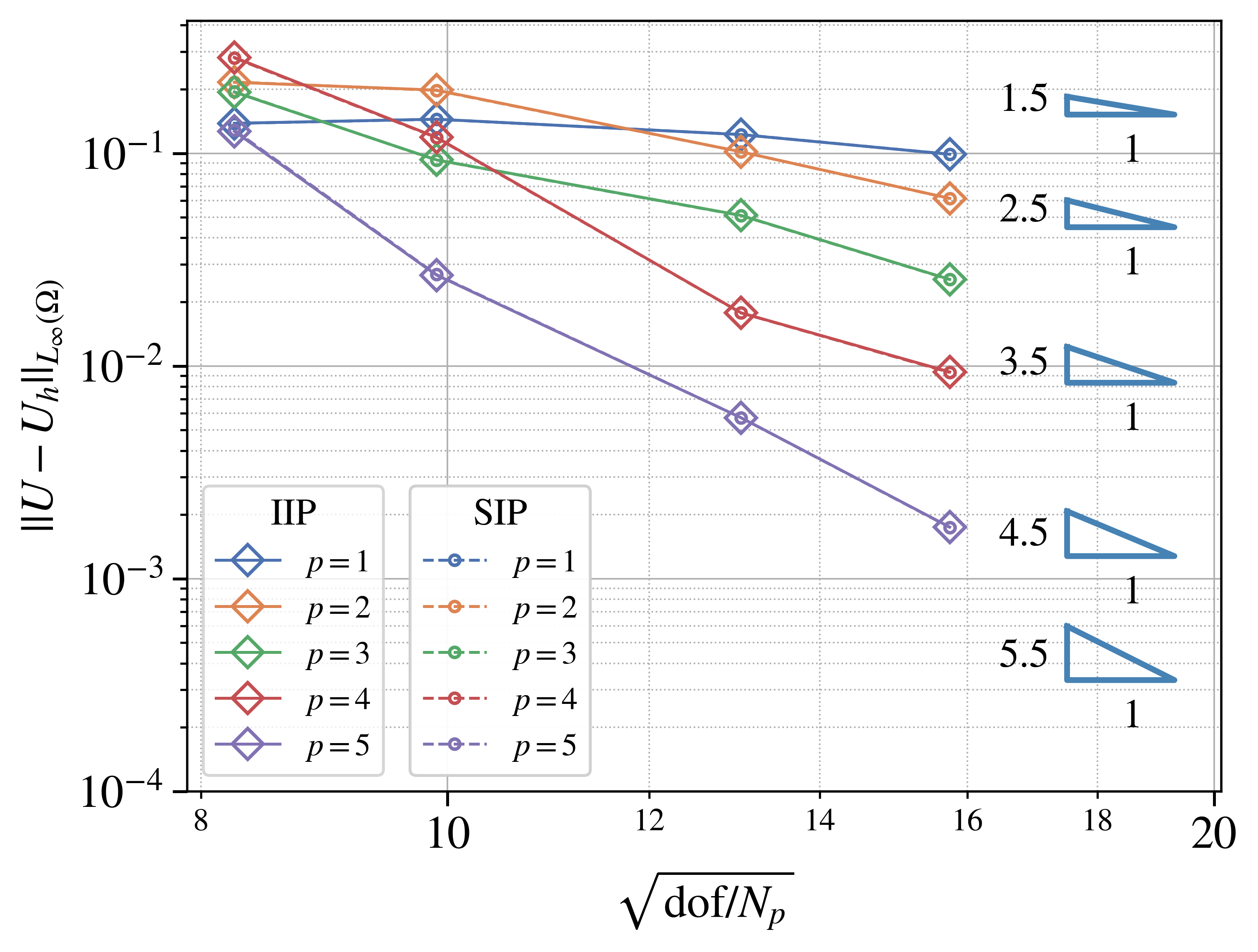}
    \end{subfigure}
    \caption{Convergence test: $L_{2}$ (top) and $L_{\infty}$ (bottom) error of $U-U_h$ 
    for $m_1$ (left) and $m_2$ (right) meshes over number of cells of the mesh.}
    \label{fig:errors_convergence}
\end{figure}
The resulting convergence rates $r_2$ and $r_{\infty}$ for $L_2(\Omega)$ and $L_{\infty}(\Omega)$ norms are provided in Tables \ref{tab:convergence_order_m1} and \ref{tab:convergence_order_m2}, respectively, together with maximum IP values.
\begin{table}
    \centering
    \begin{tabular}{c|ccc|ccc}
        \toprule 
        & 
        \multicolumn{3}{c|}{Incomplete IP} & \multicolumn{3}{c}{Symmetric IP} \\
        \midrule
        $p$ & $r_2$ & $r_{\infty}$ & $\nrm{\sigma_e}_{L_{\infty}(\Gamma)} \cdot 10^{-3}$ 
            & $r_2$ & $r_{\infty}$ & $\nrm{\sigma_e}_{L_{\infty}(\Gamma)} \cdot 10^{-3}$ \\
        \midrule 
        $1$ & $0.903$ & $0.197$ & $3.1 $ &$0.901$ & $0.198$ & $15.6 $ \\
        $2$ & $2.408$ & $0.941$ & $6.3 $ &$2.405$ & $0.938$ & $31.3 $ \\
        $3$ & $3.161$ & $1.557$ & $10.4$ &$3.162$ & $1.560$ & $51.9 $ \\
        $4$ & $4.335$ & $2.056$ & $15.6$ &$4.328$ & $2.048$ & $77.9 $ \\
        $5$ & $5.158$ & $2.611$ & $21.8$ &$5.160$ & $2.615$ & $109.0$ \\
        \bottomrule
    \end{tabular}
    \caption{Convergence test: Order of space convergence $r$ using linear regression and maximum value of the interior-penalty for $m_1$.}
    \label{tab:convergence_order_m1}
\end{table}
\begin{table}
    \centering
    \begin{tabular}{c|ccc|ccc}
        \toprule 
        & 
        \multicolumn{3}{c|}{Incomplete IP} & \multicolumn{3}{c}{Symmetric IP} \\
        \midrule
        $p$ & $r_2$ & $r_{\infty}$ & $\nrm{\sigma_e}_{L_{\infty}(\Gamma)} \cdot 10^{-3}$ 
            & $r_2$ & $r_{\infty}$ & $\nrm{\sigma_e}_{L_{\infty}(\Gamma)} \cdot 10^{-3}$ \\
        \midrule
        $1$ &$0.924$ & $0.267$ & $16.3 $ & $0.923$ & $0.266$ &$81.5 $ \\
        $2$ &$2.590$ & $1.009$ & $32.6 $ & $2.594$ & $1.012$ &$162.8$ \\
        $3$ &$4.064$ & $1.494$ & $54.2 $ & $4.060$ & $1.498$ &$271.2$ \\
        $4$ &$6.190$ & $2.757$ & $81.4 $ & $6.196$ & $2.756$ &$406.8$ \\
        $5$ &$6.749$ & $3.232$ & $113.9$ & $5.861$ & $2.745$ &$569.5$ \\
        \bottomrule
    \end{tabular}
    \caption{Convergence test: Order of space convergence $r$ using linear regression and maximum value of the interior-penalty for $m_2$.}
    \label{tab:convergence_order_m2}
\end{table}

\autoref{fig:errors_convergence} highlights that both SIP and IIP nodal RKDG approximations converge in a similar fashion without visible differences for $m_1$ and $m_2$. 
The framework provides high order accuracy $r_2(p) > p\ge 2$. Also, $r_{\infty}(p) > 2.5$ for $p\ge 4$ for $m_2$, indicating that the nodal IP RKDG retains its high order  on highly skewed meshes.

Finally, we observe that the maximum penalty values on the entire skeleton $\Gamma$, as a function of $p$, are \textit{roughly} proportional to $N_p$.

The SIP variant demands larger penalty values for stability compared to IIP, thus imposing a shorter time-step; cf. \eqref{eq:cfl_general}. Since completely analogous convergence properties have been observed for both variants for the error measures studied (\autoref{fig:errors_convergence}), for the remainder of this work, we employ the IIP variant.

\subsection{Laminar flow around stationary circular cylinder}

\label{subsection:von_karman_street}

In this section we investigate a non-vibrating circular cylinder in uniform inflow. This is a widely used validation case whereby separation occurs and a von-K\'{a}rm\'{a}n vortex street is formed downstream of the cylinder. Consequently, the shedding vortices result in highly unsteady aerodynamic loading. A review and examination of the flow past a cylinder can be found in \cite{norberg2003fluctuating}.

Initially, for this test case, a sensitivity study is needed that verifies the existence of relative convergence.
Typically, this study includes space and time refinement. However the CFL \eqref{eq:cfl_general} forces a very small time step (close to $10^{-4}\ s$ as we shall see below). For a shedding frequency $f=0.2\, Hz$ (corresponding to a Strouhal number $\mathrm{Str}:=f D / \nrm{\vect{v}_{\infty}}_2 =0.2$) these are approximately $50,000$ time-steps per period. Hence, only spatial refinement is performed and admissible time-step values are only reported.

Next, the validation with experimental and numerical results is carried out. For the validation, we employ the mesh and the polynomial degree that yield a solution insensitive to further refinement, as determined by the sensitivity study.

\subsubsection{Sensitivity study}

\begin{figure}
	\centering
	\begin{subfigure}[b]{0.61\textwidth}
		\includegraphics[width=0.95\textwidth]{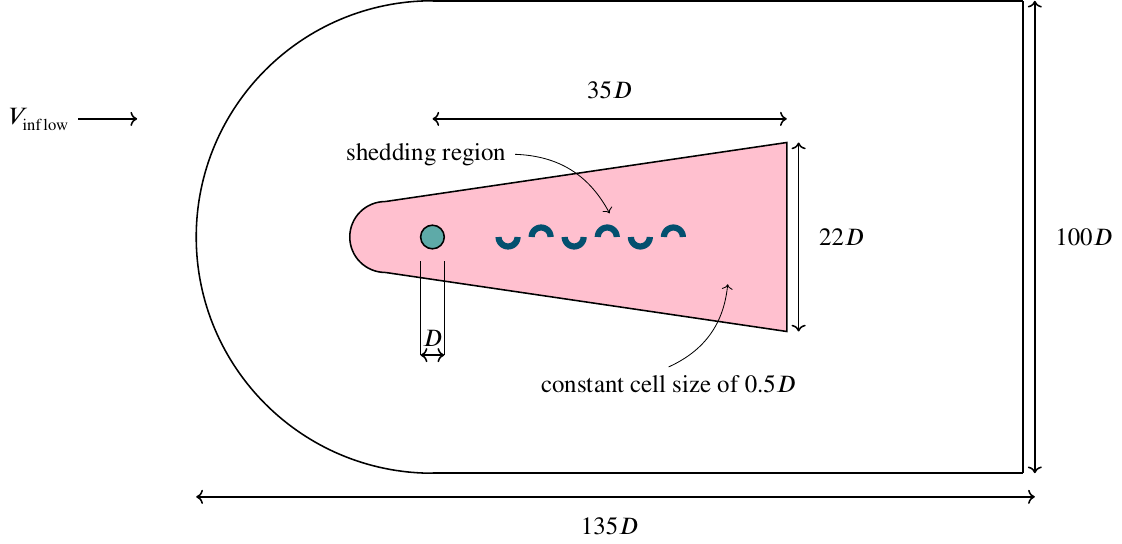}
	\end{subfigure}
	\begin{subfigure}[b]{0.38\textwidth}
		\includegraphics[width=0.98\textwidth]{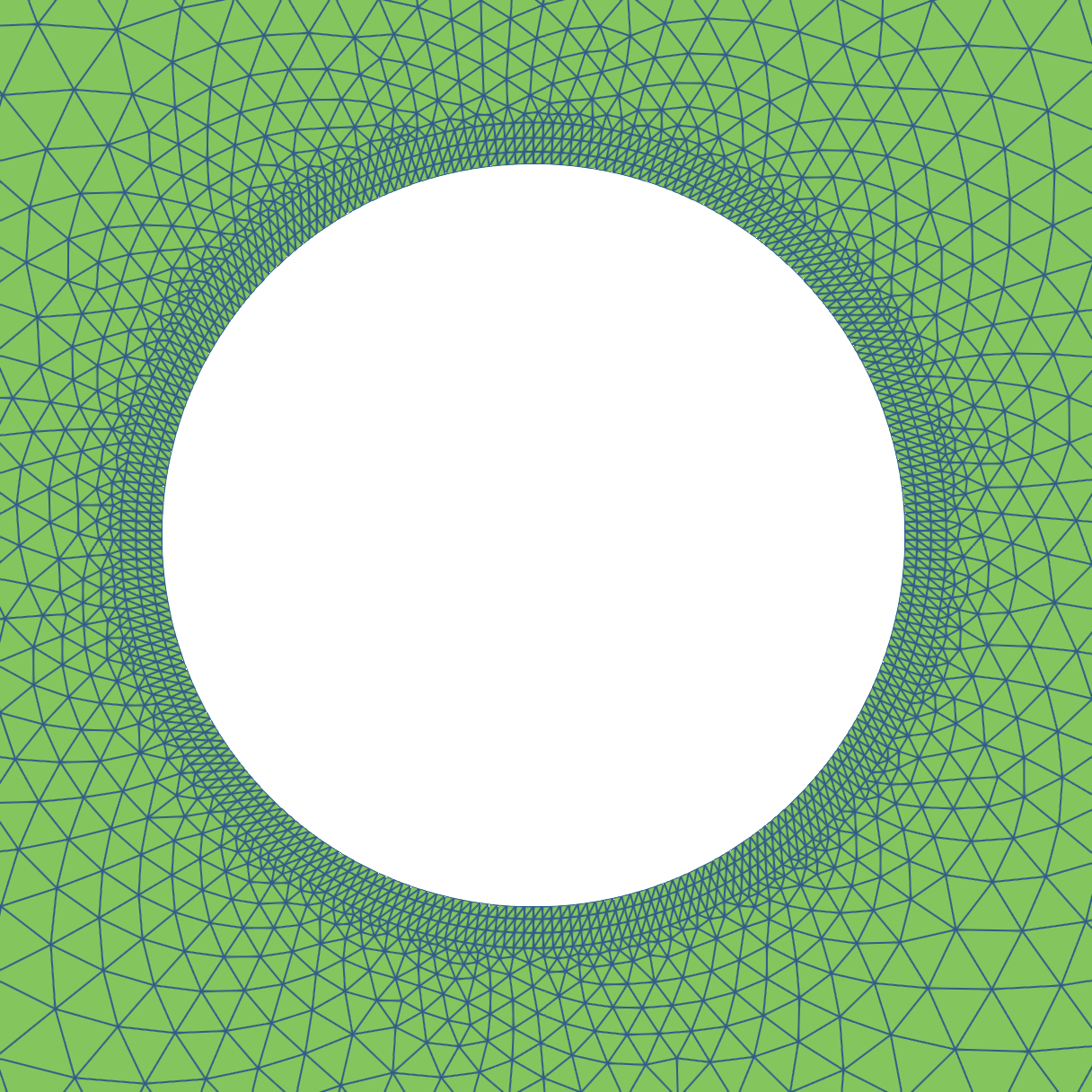}
	\end{subfigure}
	\caption{Vortex shedding case: (Left) schematic representation of $\Omega$. (Right) the mesh close to $\Gamma_{\mathrm{w}}$.}
	\label{fig:schematic_von_karman}
\end{figure}

\newcommand{\myTitle}{Sensitivity study for the vortex shedding case}
The setup for this comparative study is as follows. For a non-confined cylinder of infinite span (2D)
and diameter $D=1\, m$, the far-field region of the domain is located at $35D$ upstream, $100D$ downstream and $50D$ at each side as shown in the schematic of \autoref{fig:schematic_von_karman}.
Far-field conditions are selected to correspond to a Reynolds number $\textrm{Re}=200$, and Mach number $\textrm{Ma}=0.1$.
The mesh is constructed to contain a quasiuniform cell-size region as shown in the schematic of \autoref{fig:schematic_von_karman} to allow for the examination of numerical and physical dissipation of vortex shedding taking place there. The mesh is completed by a finer layered mesh region in the vicinity of the wall boundary to capture the boundary layer, as shown in \autoref{fig:schematic_von_karman} (right panel), and a far-field region with adjacent cells of size $h_K=8\, m$. Further details are listed in \autoref{tab:mesh_particulars}.
\begin{table}
	\centering
	\begin{tabular}{lr}
		\toprule
		Description                            & value       \\
		\midrule
		Number of cells                        & $20,039   $ \\
		cell length at wall                    & $12\, mm  $ \\
		height of first cell from wall         & $17\, mm  $ \\
		total height of layered mesh from wall & $56\, mm  $ \\
		cell length at shedding region         & $500\, mm $ \\
		cell length at far-field               & $8-8.5\, m$ \\
		\bottomrule
	\end{tabular}
	\caption{Vortex shedding case: mesh characteristics.}
	\label{tab:mesh_particulars}
\end{table}

The sensitivity study begins by increasing $p$, for $1\le p \le 5$.
A convergent behaviour on the relative difference $\delta_p \mathrm{f} := \mathrm{f}(\cdot,p+1) - \mathrm{f}(\cdot,p)$ of a result $\mathrm{f}$ is achieved if $\delta_p \mathrm{f}$ is decreasing with respect to $p$.

For assessment, we use the vorticity field, which is usually indicative of diffusion when observed across the downstream direction. The examination is both visual and also quantitative.
\begin{figure}
	\centering
	\begin{subfigure}[b]{0.99\textwidth}
		\includegraphics[width=\textwidth]{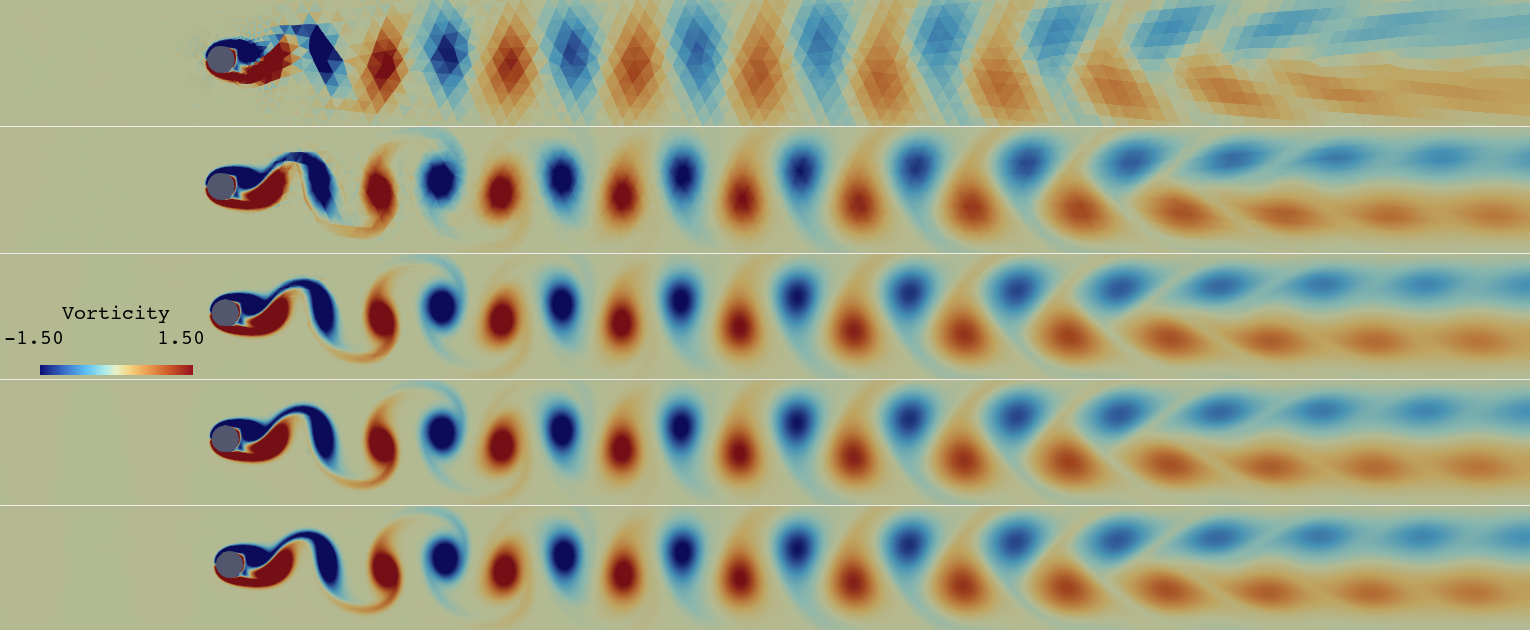}
	\end{subfigure}
	\begin{subfigure}[b]{0.99\textwidth}
		\includegraphics[width=\textwidth]{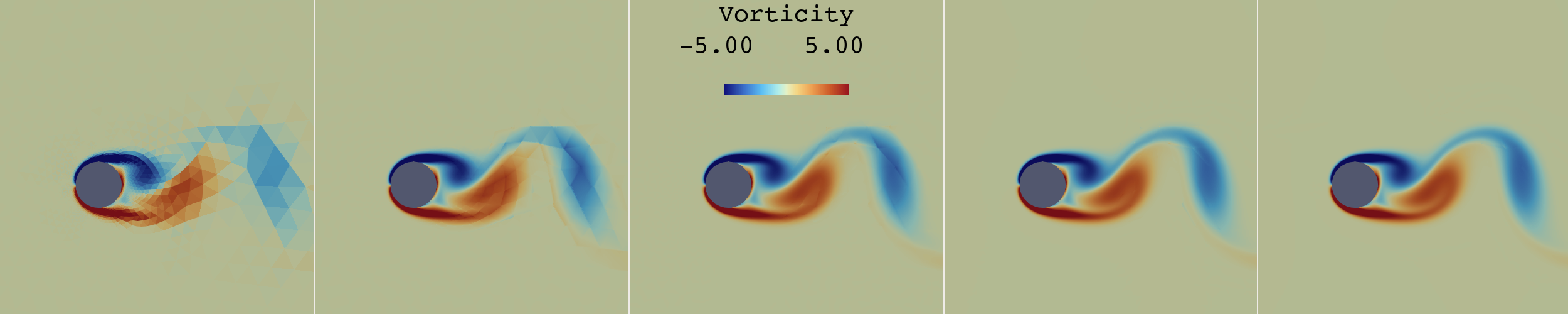}
	\end{subfigure}
	\caption{\myTitle: Vorticity contour plots at $\mathrm{Re}=200$.
	}
	\label{fig:vort_contour_sens}
\end{figure}

In \autoref{fig:vort_contour_sens} a qualitative comparison of results for $p=1,2,3,4,5$ is presented. These snapshots refer to $t> 190\ s$, past which the behaviour is periodic, and the results are synchronised in shedding phase for the purpose of comparison. From the vorticity field in the downstream direction, it is evident that the phenomenon appears to be spatially resolved for $p\ge 3$.

\begin{figure}
	\centering
	\begin{subfigure}[b]{0.92\textwidth}
		\centering
		\includegraphics[width=0.98\textwidth]{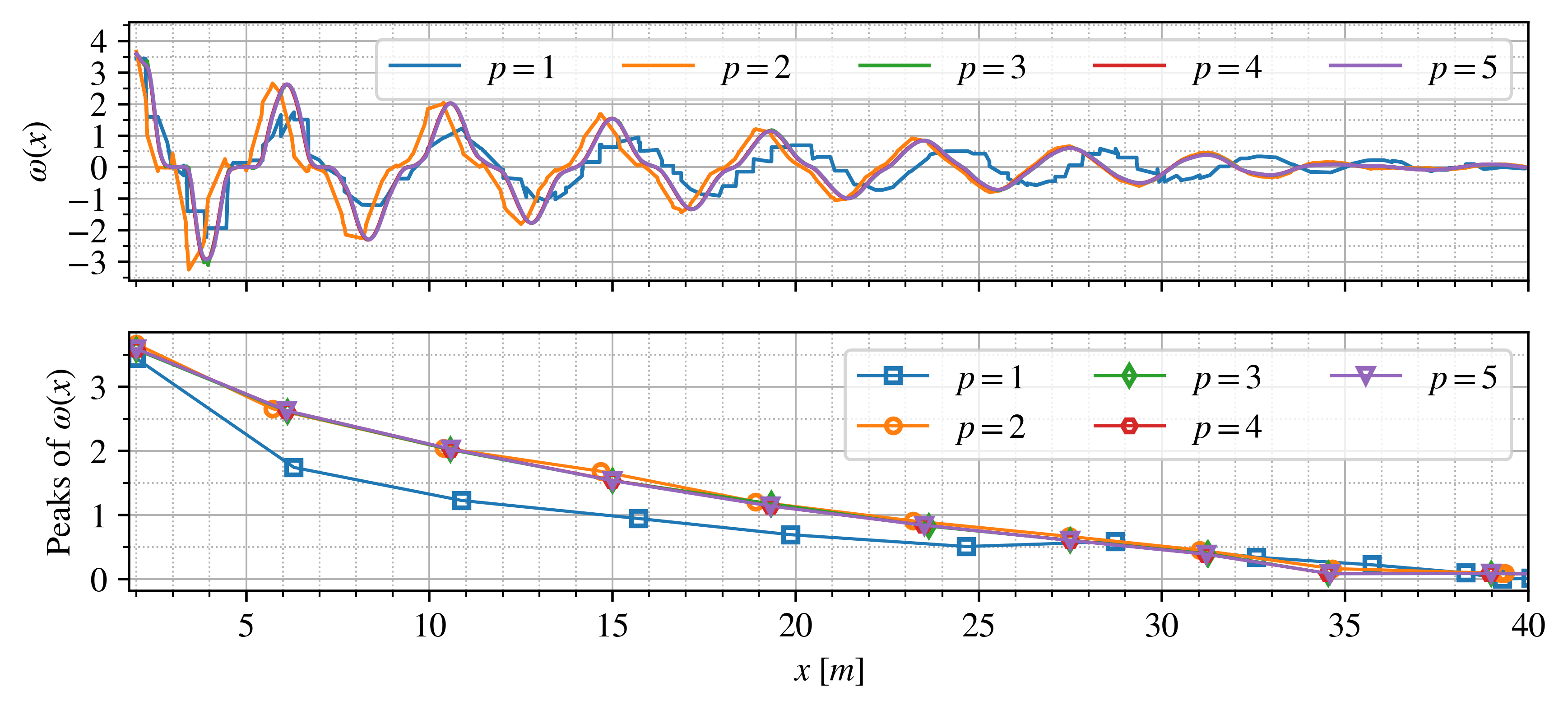}
	\end{subfigure}
	\begin{subfigure}[b]{0.90\textwidth}
		\centering
		\includegraphics[width=0.98\textwidth]{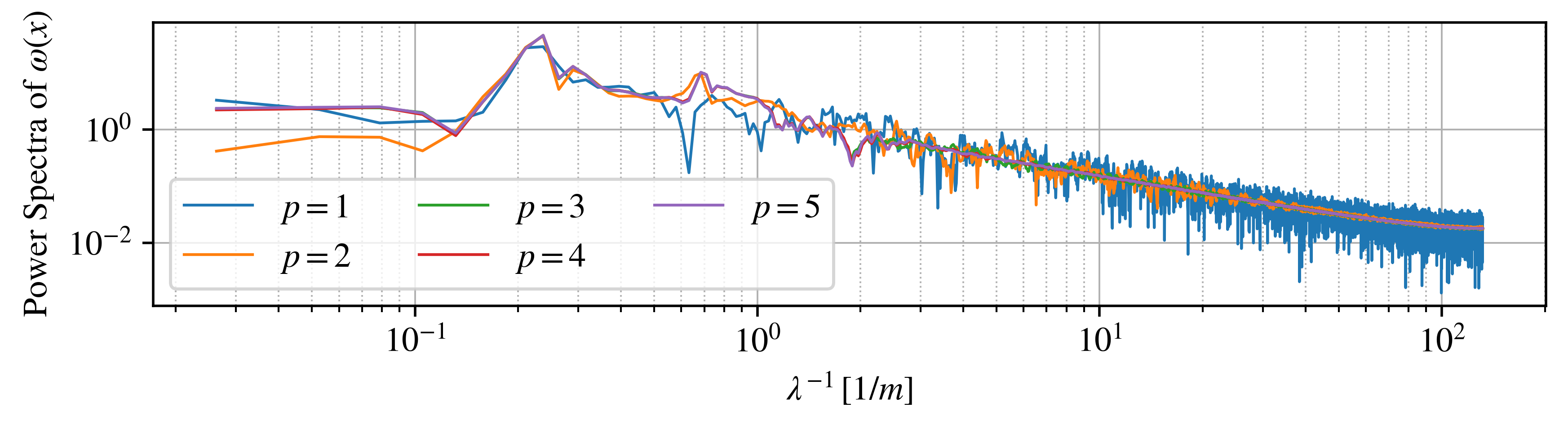}
	\end{subfigure}
	\caption{\myTitle: (Top) Vorticity $\omega(x)$ for $x\in I_x$. (Middle) Peaks of $\omega(x)$. (Bottom) Magnitude of the (space) discrete Fourier transform of $\omega(x)$.}
	\label{fig:vorticity_peaks_and_dft}
\end{figure}

Quantitatively, vorticity $\omega(x)$ is extracted on a line segment across the downstream vector $[(40,0) - (2,0)]^T\, m$ and the relative error $\delta_p \omega$ is measured. In \autoref{fig:vorticity_peaks_and_dft}, we display the vorticity peaks and the power spectra's magnitude of this vorticity from a discrete Fourier transform. We observe visible differences at the vorticity's amplitude and the wavelength between two peaks for $p=1$ and $p=2$, which vanish for $p> 2$. This is also observed in \autoref{tab:vorticity_norms} numerically, whereby values of $\nrm{\delta_p\omega}_{L_2\pth{I_x}}$, $\nrm{\delta_p \omega}_{L_{\infty}\pth{I_x}}$ with $I_x=(2,40)$, for $p\ge 3$ are at least a magnitude smaller than those from $p=1$ or $p=2$; a remark that validates the visual finding of \autoref{fig:vort_contour_sens}.
\begin{table}
	\centering
	\begin{tabular}{rll}
		\toprule
		$p$ & $\nrm{\delta_p\omega}_{I_x}$ & $\nrm{\delta_p\omega}_{L_{\infty}\pth{I_x}}$ \\
		\midrule
		$1$ & $0.683$                      & $2.194$                                      \\
		$2$ & $0.517$                      & $3.056$                                      \\
		$3$ & $0.025$                      & $0.188$                                      \\
		$4$ & $0.030$                      & $0.169$                                      \\
		\bottomrule
	\end{tabular}
	\caption{\myTitle: $p$-relative errors of stream-wise vorticity.}
	\label{tab:vorticity_norms}
\end{table}

Next, sectional lift, $C_L$, and drag, $C_D$, coefficient differences are presented, whereby
\begin{align}
	 & C_L := \frac{\mathrm{Lift}}{\frac{1}{2}\rho \vect{v}^2_{\infty} D},
	 & C_D := \frac{\mathrm{Drag}}{\frac{1}{2}\rho \vect{v}^2_{\infty} D},
	\label{eq:lift_drag}
\end{align}
with the lift and drag forces considered to be perpendicular and parallel to the inflow vector $\vect{v}_{\infty}$, respectively.
The average ``$\overline{\cdot}$'' and root mean squared (RMS) values of $C_L$, $C_D$ are evaluated ahead of the transient state, i.e., for $150\le t\le 200\, s$, and the results are provided in \autoref{tab:mean_rms_cl_cd}. \autoref{tab:mean_rms_cl_cd} suggests that, with the exception of $\overline{C}_L$ which is insensitive for $p=4,5$, all other quantities have a converging trend with increasing $p$.
\begin{table}
	\centering
	\begin{tabular}{rrrrr}
		\toprule
		$ (\ \delta_p\overline{C}_L$ & $\delta_p C_{L,\mathrm{RMS}}$ & $\delta_p\overline{C}_D$ & $\ \delta_p C_{D,\mathrm{RMS}}\ )$ & $\times10^{-3}$ \\
		\midrule
		$-3.882$                     & $3.3297$                      & $19.5921$                & $19.5884$                          &                 \\
		$-0.826$                     & $1.5710$                      & $ 0.7368$                & $ 0.7411$                          &                 \\
		$-0.682$                     & $0.2163$                      & $ 0.1801$                & $ 0.1806$                          &                 \\
		$-0.785$                     & $0.0015$                      & $-0.0081$                & $-0.0083$                          &                 \\
		\bottomrule
	\end{tabular}
	\caption{\myTitle: differences for $C_D$ and $C_L$ during $p-$sensitivity.}
	\label{tab:mean_rms_cl_cd}
\end{table}

From the above results, we conclude that at $\mathrm{Re}=200$, sufficient flow resolution can be achieved with $p=3$. In the following validation $p=3$ is used.

Further, we report the maximum permissible time-step values. These were found heuristically and the results are tabulated in \autoref{tab:dt_max}. \autoref{tab:dt_max} suggests that the system becomes computationally demanding for $p \ge 5$ in terms of numbers of iterations and one may have to resort to more sophisticated strategies.
\begin{table}
	\centering
	\begin{tabular}{lr}
		\toprule
		$p$ & $\delta t \cdot 10^{-4} $ \\
		\midrule
		$1$ & $3.766$                   \\
		$2$ & $1.815$                   \\
		$3$ & $0.917$                   \\
		$4$ & $0.427$                   \\
		$5$ & $0.308$                   \\
		\bottomrule
	\end{tabular}
	\caption{\myTitle: maximum permissible time-step values for up to $p=9$ polynomial order.}
	\label{tab:dt_max}
\end{table}

\subsubsection{Validation}
\label{subsubsection:validation_von_karman}
We now validate against existing results from the literature. The experimental data refer to \cite{norberg1994experimental} and the numerical data are collected from various studies. For $100\le \mathrm{Re}\le 250$, the Strouhal number and the (mean) drag coefficient are used as comparison measures.

We note that both experimental and numerical data refer to incompressible conditions. Selecting $\mathrm{Ma}=0.1$ on $\Gamma_{\mathrm{far}}$, we eliminate compressibility effects.
\begin{figure}
	\centering
	\begin{subfigure}[b]{0.49\textwidth}
		\centering
		\includegraphics[width=0.98\textwidth]{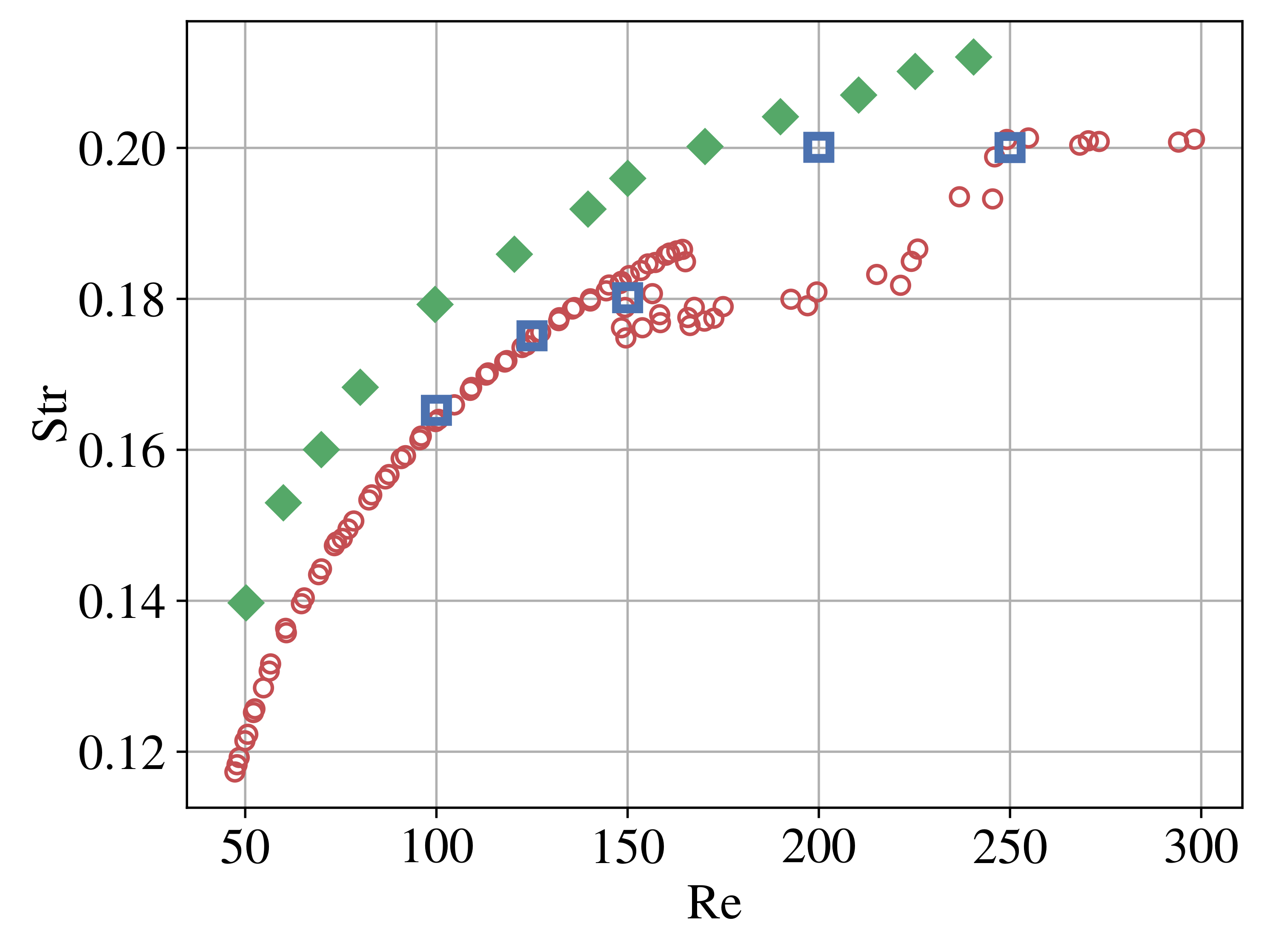}
	\end{subfigure}
	\begin{subfigure}[b]{0.49\textwidth}
		\includegraphics[width=0.98\linewidth]{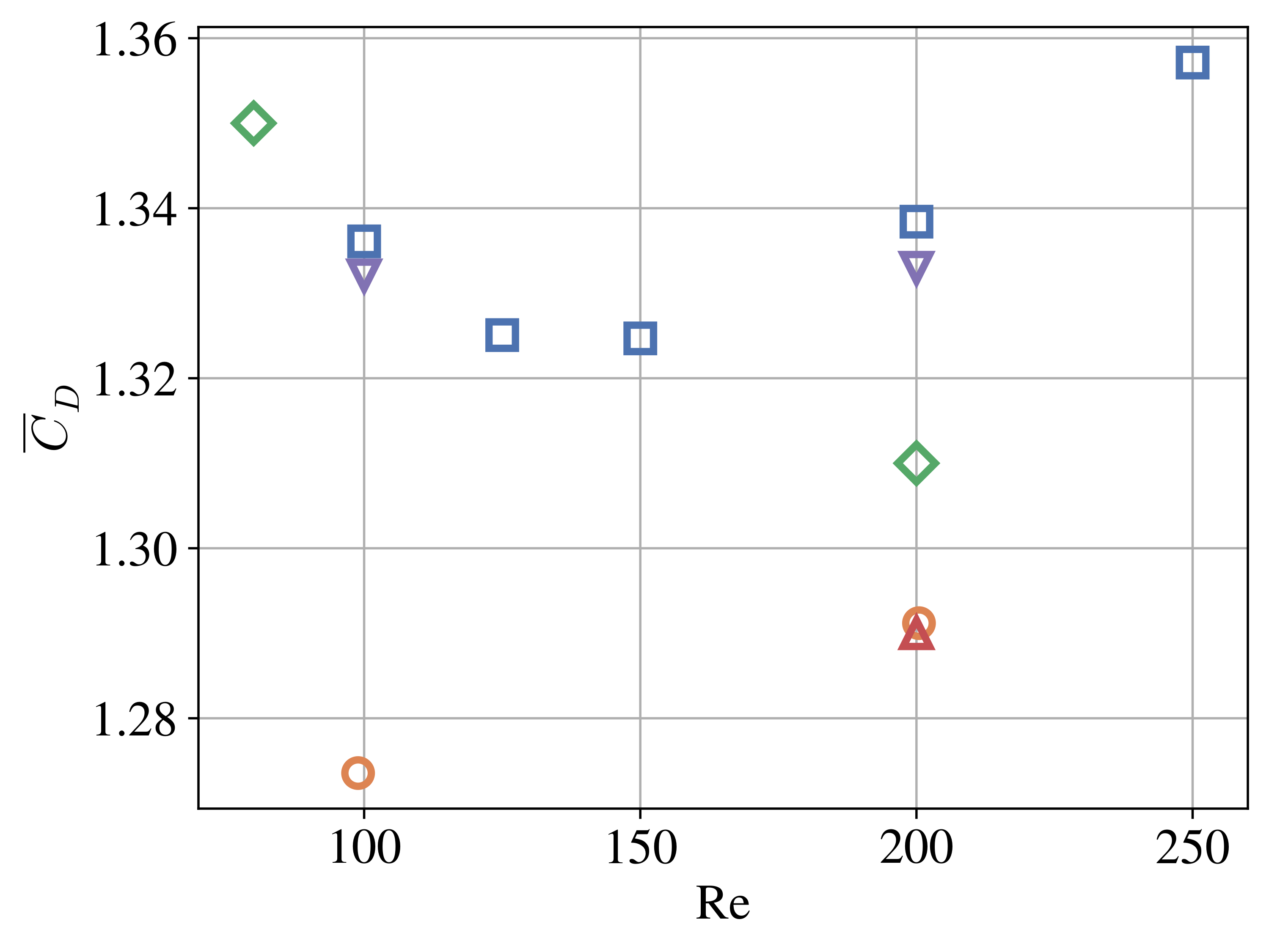}
	\end{subfigure}
	\caption{Validation of vortex shedding case: (Left) Strouhal number over Reynolds number: ($\square$) present study; ($\circ$) exp. from \cite{norberg1994experimental}; ($\blacklozenge$) num. from \cite{karniadakis1989frequency}. (Right) drag coefficient over Reynolds number: ($\square$) present study; ($\circ$) num. from \cite{braza1986}; ($\blacklozenge$) num. from \cite{Franke1990237}; ($\bigtriangleup$) num. from \cite{Rogers1991}; ($\bigtriangledown$) num. from \cite{Posdziech2007479}.}
	\label{fig:shedding_validation}
\end{figure}

In \autoref{fig:shedding_validation}, the measured Strouhal number is in agreement with the experiment data, especially for $\mathrm{Re}< 200$. For $\mathrm{Re}=200$ there are visible differences that are attributed to the effects of three-dimensional interactions. These are significant for $Re \ge 190$; as stated in \cite{norberg1994experimental} and explained thoroughly in the review \cite{ForouziFeshalami2022}. Further, the values of $\overline{C}_D$ for $100\le \mathrm{Re} \le 200$ are also in agreement with other numerical studies from the literature, indicating the validity of our results for the vortex shedding past a circular stationary cylinder.

\subsection{Elastically mounted cylinder}

\label{subsection:elastically_mounted_cylinder}

Building upon the validation results for stationary cylinders presented above, we now extend our investigation to a more complex fluid-structure interaction problem. More specifically, a circular cylinder is now elastically mounted on a spring-damper system. 
In this coupled system, vortices shed by the solid body cause movement around the oscillation's equilibrium, and this motion then modifies the ongoing vortex shedding process. This concept problem falls into the wider category of vortex induced vibration (VIV) problems. 
Here, a single body's movement is recreated by a rigid mesh motion, allowing for the present numerical framework to address this case also.

\subsubsection{Definition of the problem}
\begin{figure}
    \centering
    \includegraphics[width=0.75\textwidth]{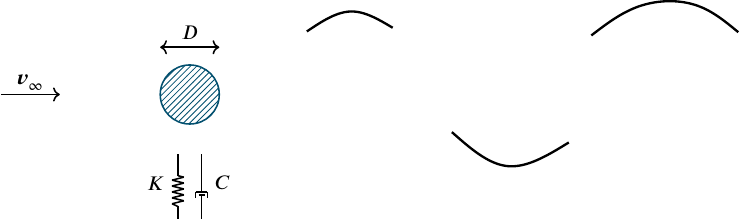}
    \caption{Schematic representation of a 1-dof elastically-mounted oscillating cylinder.}
    \label{fig:schematic_elastically_mounted}
\end{figure}

A cylinder is modelled as a one degree of freedom (\textit{dof}) oscillator, attached vertically to a spring and a damper, as seen in the schematic in \autoref{fig:schematic_elastically_mounted}. The moving dynamics equation \eqref{eq:1dof_oscillator} is used with coefficients:
\begin{gather*}
    M_{\mathrm{r}} = 1,\quad C_{\mathrm{r}} = 4\pi f_n \xi,\quad K_{\mathrm{r}} = \Big(1 + \frac{1}{m^*}\Big)(2\pi f_n)^2,\quad F_{\mathrm{r}} = \frac{\text{Lift}}{m};
\end{gather*}
the remaining symbols are defined in \autoref{tab:1d_oscillator}.
\begin{table}[ht]
    \centering
    \ra{1.1}
    \begin{tabular}[c]{lll}
        \toprule
        Definition                                               & description                            \\
        \midrule        
        $m$                                                      & mass                                   \\
        $m_f := \rho \pi D^2/4$                                  & mass of displaced volume               \\
        $m^* := m/m_f = 1$                                       & reduced mass                           \\
        $\mathrm{k}$                                             & spring coefficient                     \\
        $\xi=0.01$                                               & damping coefficient                    \\
        $f_n := \frac{1}{2\pi} \sqrt{\frac{\mathrm{k}}{m+m_f}} $ & natural frequency including added mass \\
        $U^* := \frac{U_{\infty}}{f_n D}$                        & reduced velocity                       \\
        \bottomrule
    \end{tabular}
    \caption{Elastically-mounted cylinder: Description of parameters.}
    \label{tab:1d_oscillator}
\end{table}

Our aim is to study the resulting response by (parametrically) changing the spring and damping strength. Simulations are carried out for $U^* \in [3.0,7.0]$ using the mesh from the validation in \autoref{subsection:von_karman_street} with $p=3$ and the time-step dictated by \eqref{eq:cfl_general}. It is noted that, given a $U^*$, $f_n$ can be evaluated, which is then used to derive the spring and damping strengths via simple algebra on relations from \autoref{tab:1d_oscillator}.

According to \cite{leontini2006beginning,khalak1999motions}, by decreasing the spring and damping constants (i.e., increasing $U^{*}$,) a secondary frequency mode arises, which becomes dominant for $U^* \ge 3.6$. Then the first two frequencies create a branching region, where oscillation characteristics change abruptly. The existence of branching will be one of the key validating points of this test case.

\subsubsection{Numerical results}
To compare with available data from the literature, we first define the sampling time interval $I_a = (60,210)\ s$. This is chosen to include at least $30$ shedding cycles of the corresponding stationary cylinder problem ($\mathrm{Str} = 0.2$), excluding the transient interval $(0,60)\ s$.

Then we consider three quantities:
\begin{inparaenum}[1)]
    \item the primary frequency of oscillation $f_{\mathrm{prim}}$ calculated as in \autoref{subsection:von_karman_street} using a DFT on $I_a$;
    \item the maximum lift coefficient $\nrm{C_L}_{L_{\infty}(I_a)}$; and
    \item the maximum oscillation amplitude $\nrm{A}_{L_{\infty}(I_a)}$.
\end{inparaenum}

\begin{figure}
    \centering
    \includegraphics[width=0.99\textwidth]{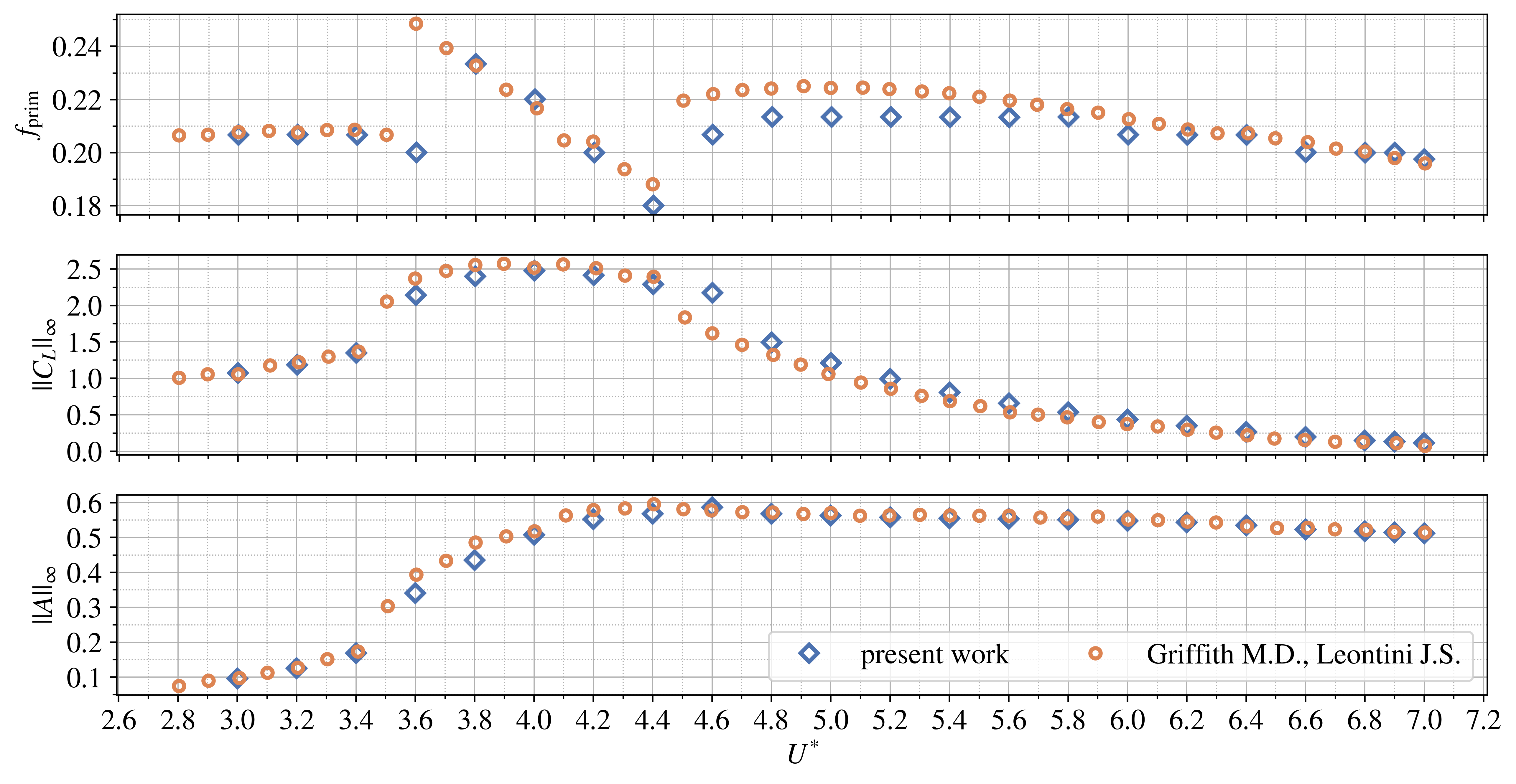}
    \caption{Elastically-mounted cylinder: comparison of the oscillation's dominant frequency (upper), the maximum value of lift coefficient (middle) and the maximum value of response's amplitude (lower figure) over the normalised natural frequency.}
    \label{fig:comparison_with_u_star}
\end{figure}
\begin{figure}
    \begin{center}
        \includegraphics[width=0.99\textwidth]{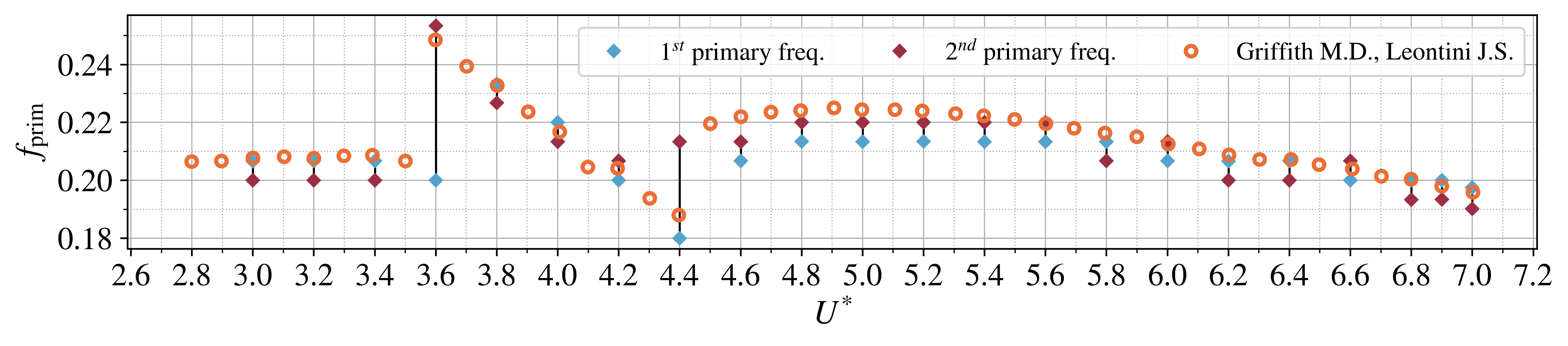}
    \end{center}
    \caption{Elastically-mounted cylinder: comparison with the first two dominant oscillation frequencies.}\label{fig:first_second_with_u_star}
\end{figure}
\begin{figure}
    \centering
    \includegraphics[width=0.90\textwidth]{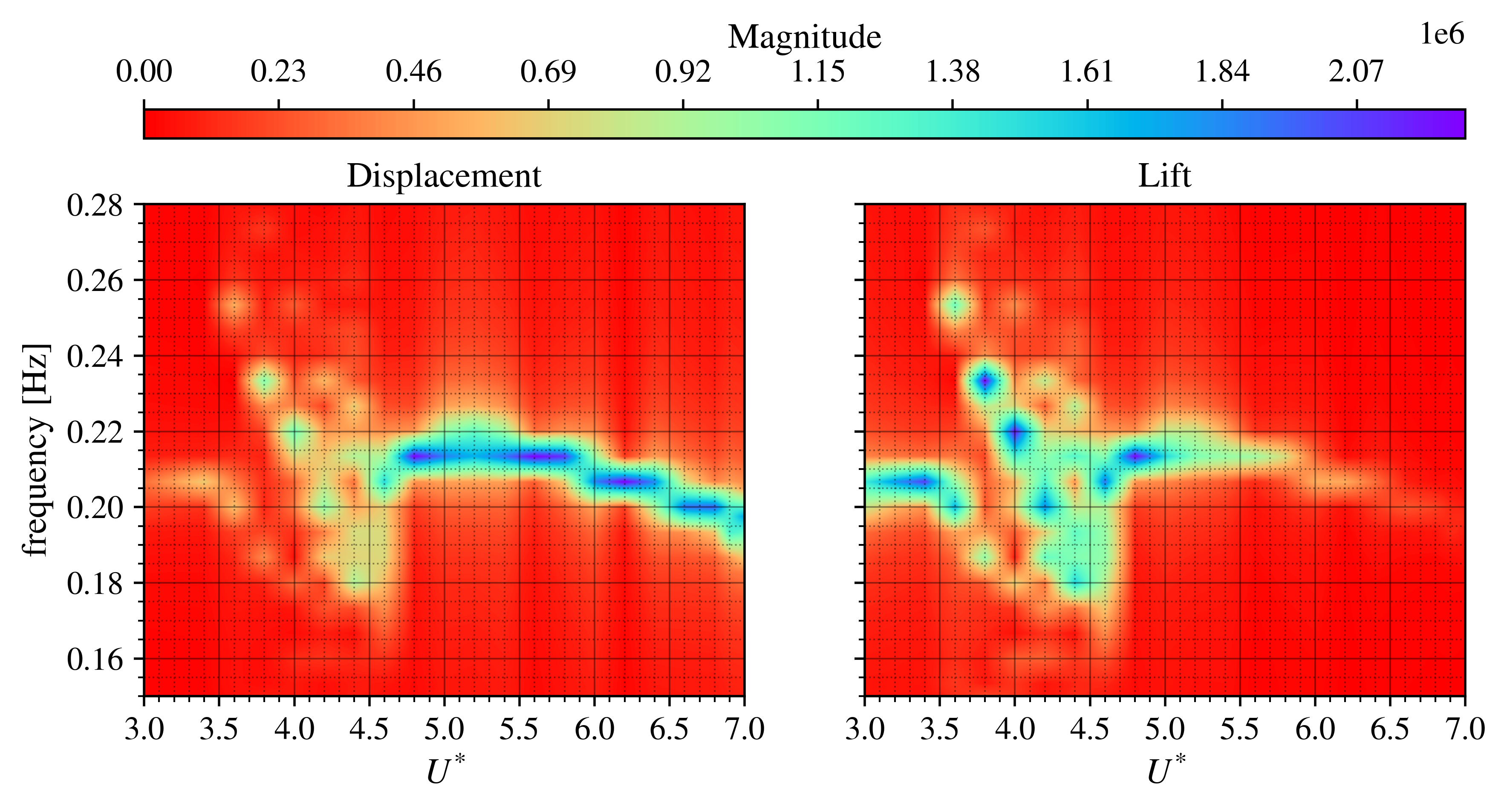}
    \caption{Elastically-mounted cylinder: contour plot of power spectra of cylinder's displacement (left) and lift force (right) for $U^*\in(3.0,7.0)$.}
    \label{fig:spectral_amp}
\end{figure}
The results from the nodal IP RKDG are compared in \autoref{fig:comparison_with_u_star} with numerical data from a spectral element code from \cite{griffith2017sharp} and presented in detail in \cite{leontini2006beginning}.

From \autoref{fig:comparison_with_u_star}, we observe that $\nrm{C_L}_{L_{\infty}(I_a)}$ and $\nrm{A}_{L_{\infty}(I_a)}$ values are in very good agreement with numerical data of \cite{griffith2017sharp}, and follow the same trend with increasing $U^*$. 
Also, the sudden change in $f_{\mathrm{prim}}$ at the start and at the end of the branching region for $U^* =(3.6,4.4)$ agrees with the region proposed by available data. However, there are two significant differences in the dominant frequency.

We believe these exist due to the problem at hand containing more than one distinct oscillation modes, and for $U^* = 3.6$ and $U^* = 4.4$ the first two largest modes of the power spectra, as seen in \autoref{fig:first_second_with_u_star}, concern very different frequency values. The magnitude of these modes is visualized in \autoref{fig:spectral_amp} both for the lateral displacement and the lateral force, as a function of $U^*$. From \autoref{fig:spectral_amp}, it is evident that the dominant oscillation frequencies refer to modes with very close magnitudes both for $U^* = 3.6$ and $U^* = 4.4$ and, thus, $f_{\mathrm{prim}}$ in \autoref{fig:comparison_with_u_star} may differ significantly from the one acquired from \cite{griffith2017sharp}, since it refers to the second largest oscillation frequency.

From the contours from \autoref{fig:spectral_amp} for $U^* > 4.5$, in which the higher-frequency mode becomes more dominant, the magnitude of spectral modes of $C_L$ (right figure) is smaller compared to this of oscillation's displacement (left figure), indicating that the shedding intensity drops, while the oscillation intensity increases. This observation justifies the increase in amplitude shown in \autoref{fig:comparison_with_u_star}.

\begin{figure}
    \centering
    \includegraphics[width=0.9\linewidth]{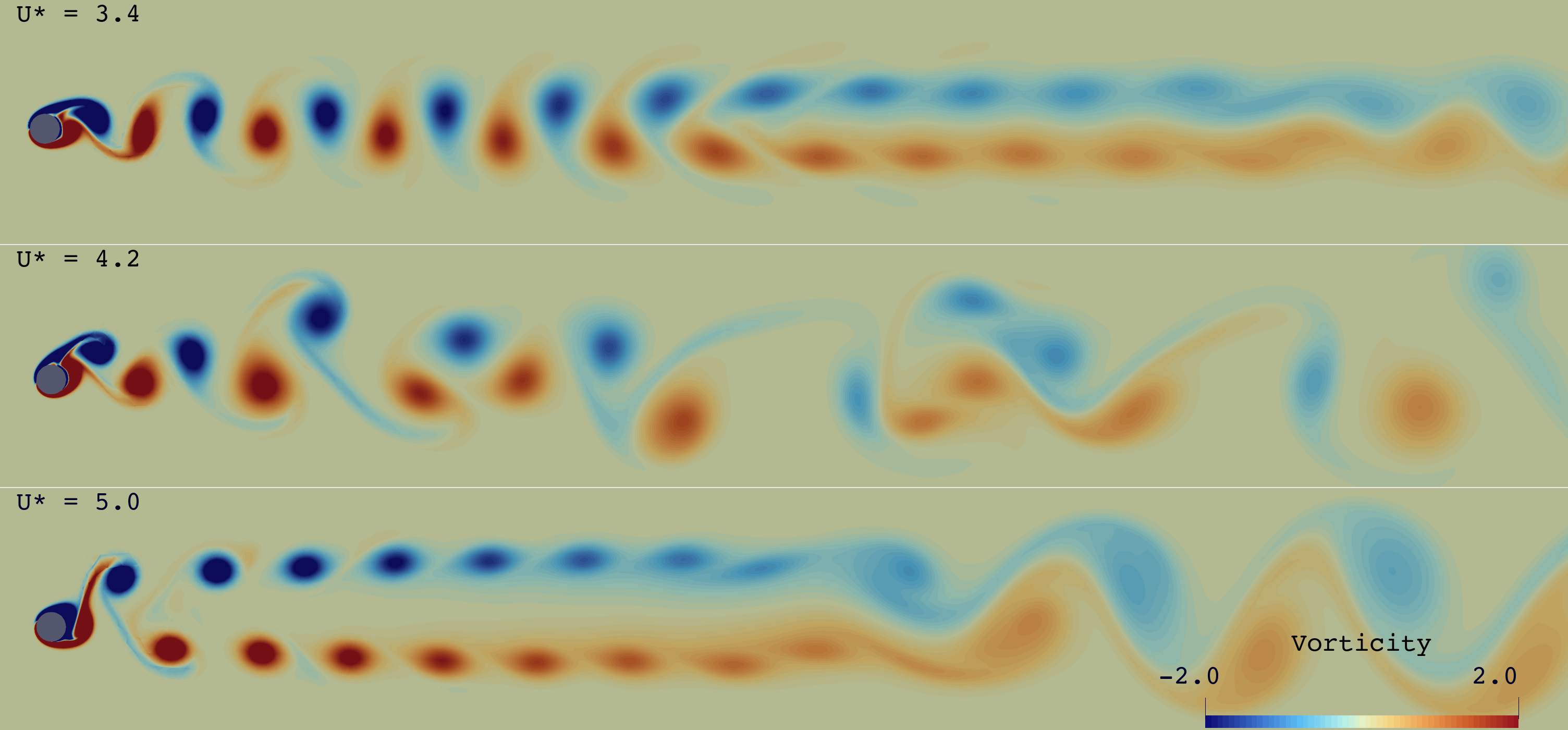}
    \caption{Elastically-mounted cylinder: vorticity contour for $U^*\in \{3.4, 4.2, 5.0\}$ over a portion of $\Omega$ at $t\approx 198\ s$.}
    \label{fig:3_u_star_vort}
\end{figure}
\begin{figure}
    \centering
    \includegraphics[width=0.9\linewidth]{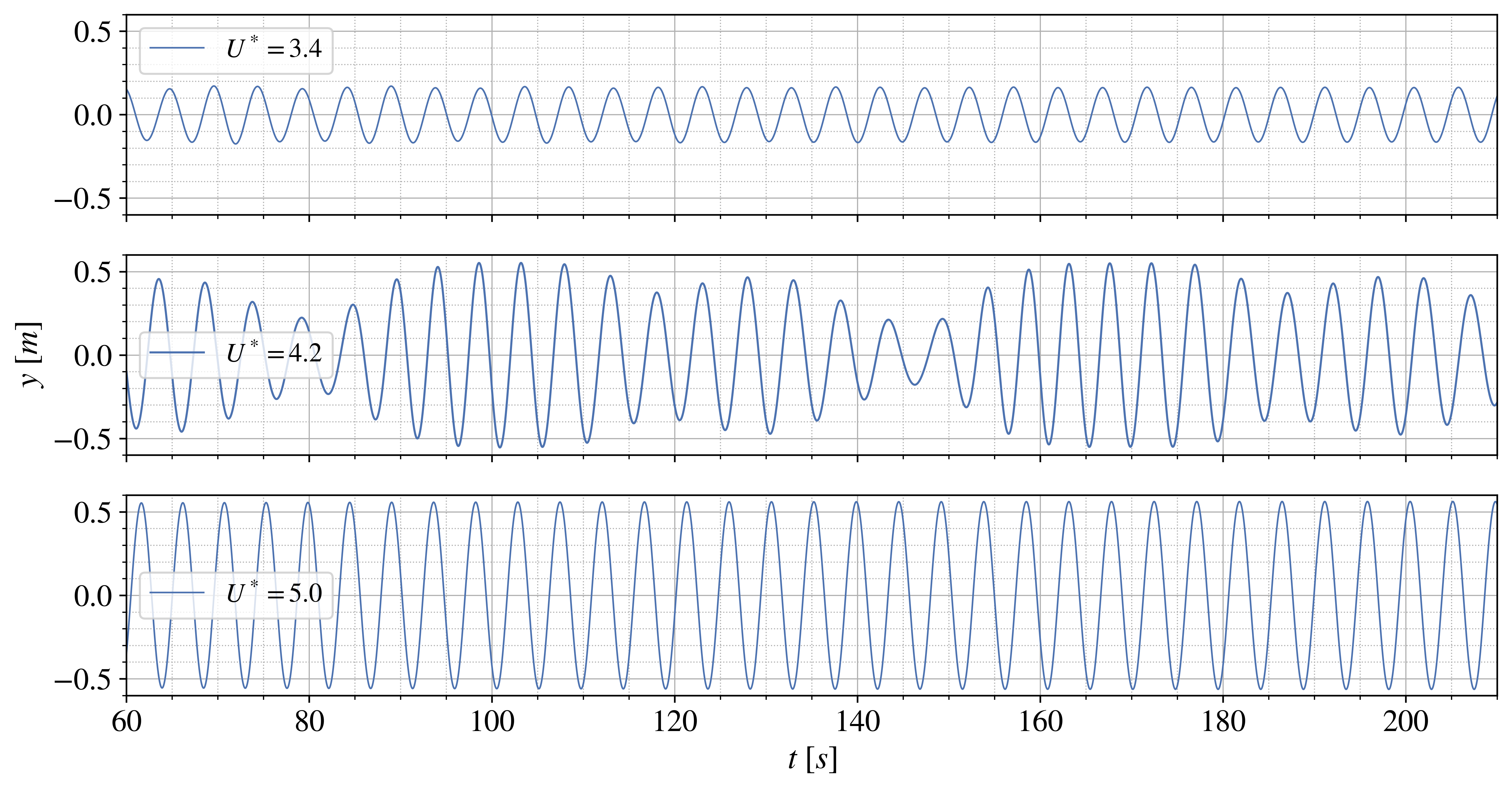}
    \caption{Elastically-mounted cylinder: response displacement over time for $U^*\in \{3.4, 4.2, 5.0\}$.}
    \label{fig:3_u_star_heave}
\end{figure}
Further, in \autoref{fig:3_u_star_vort} and \autoref{fig:3_u_star_heave}, we present a comparison among the three values of $U^*= 3.4, 4.2,5.0$. Each value refers to a different oscillation response region from $(3.0,3.6), (3.6,4.4)$ and $(4.4,7.0)$.
The contours in \autoref{fig:3_u_star_vort} reveal that the shedding type differs among these three regions. More specifically, for $U^*=3.4$, a $2S$ shedding type is observed; resembling a stationary cylinder shedding; $U^* = 4.2$ results in a $P+S$ shedding type. Lastly, for $U^*=5.0$ a $2P$ shedding type dictates the downstream image. We note that, \textit{S} and \textit{P} refer to single and paired vortices, respectively.

\section{Conclusions}

In this study, we propose and test extensively a computational framework based on a parallel implementation of a nodal Runge-Kutta interior-penalty discontinuous Galerkin method, with careful control of the CFL condition constants, which is successfully applied to unsteady test cases including a VIV problem. The methodology is shown to yield high order accuracy and it is validated against experimental and numerical data wherever available. The methodology hinges on a new, careful choice of the interior penalty parameter for the non-linear diffusion, which, together with the nodal basis implementation, the explicit RK time discretisation and the coupling with rigid dynamics, lead to an efficient solver. The solver is shown to seamlessly scale up to massively parallel rank numbers. The approach extends verbatim to three spatial dimensions, therefore, paving the way for exploring 3D CFD cases providing high fidelity results in the future.

\section{Acknowledgements}
The authors acknowledge \textit{EuroHPC Joint Undertaking} for awarding the
project with ID:EHPC-DEV-2024D03-032,
access to MeluXina at LuxProvide, Luxembourg.
This facility was used for all the test cases of this study.
Also, the authors acknowledge BETA CAE Systems for providing the mesh generation
software: ANSA, which was used for the entire meshing.
All visual post processing results, including flow field
contours, have been created using Paraview, a free and open-source software.
\section{Funding}
The research work was supported by the Hellenic Foundation for Research and Innovation (HFRI) under the 5th Call for HFRI PhD Fellowships (Fellowship Number: 20716), which involves the first author. The second author gratefully acknowledges the financial support of EPSRC (grant number EP/W005840/2).

\appendix

\section{Appendix: On the CFL condition}\label{CFL_Heuristics}
We now illustrate the heuristic choice of the CFL condition used in this work, for a corresponding linearized scalar problem, focusing on its dependence on the problem coefficients and on the penalty parameter.  More specifically, we seek $u:I\times \Omega\subset \mathbb{R}^d\to\mathbb{R}$ solving
\begin{equation}\label{PDE_scalar}
	u_t-\nabla\cdot ( \mathcal{G}\nabla u -\mathbf{b} u)	=0,
\end{equation}
on $I\times \Omega$ with no-slip boundary conditions for simplicity, $\mathcal{G}$ diffusion tensor,  and advection field $\mathbf{b}$, independent of $u$.

For this section \emph{only}, on a face $e=\partial K_1\cap\partial K_2$ shared by two elements $K_1,K_2\in\mathcal{T}$, we denote by $\jumpsc{w}|_e:=(w|_{K_1}\mathbf{n}_1+w|_{K_2}\mathbf{n}_2)|_e$ the jump across $e$ with $\mathbf{n}_i$ denoting the unit outward normal of $K_i$ to $e$; when $e\subset \partial K\cap\Gamma_b$, we set  $\jumpsc{w}|_e:=(w|_{K}\mathbf{n})|_e$ with $\mathbf{n}$ denoting the unit outward normal to $\Gamma_b$. Further, on each face $e$, let $\lfloor w\rfloor|_e:=\lim_{\epsilon\to 0^+}(w(\cdot+\epsilon\mathbf{b})-w(\cdot-\epsilon\mathbf{b}))|_e$ denote the upwind jump, which may differ from $\jumpsc{w}$ up to a sign.
We denote by $(\cdot,\cdot)_\omega$ and $\|\cdot\|_\omega$, for $\omega\subset \Omega$, the $L^2(\omega)$ inner product and norm, respectively, integrating over the corresponding measure of $\omega$.

The space discretization by the IPDG method, then, reads: find $u_h\in V_h$, such that
\begin{equation}\label{space_linear_DG}
	\int_\Omega (u_h)_t v_h\dx+B(u_h,v_h)=0,
\end{equation}
for all $v_h\in V_h$, where $B(u_h,v_h):=B_{\rm d}(u_h,v_h)+B_{\rm a}(u_h,v_h)$, for $\theta\in [-1,1]$, with
\[
	B_{\rm d}(w,v):=\int_{\Omega}\mathcal{G}\nabla_h w\cdot \nabla_h v \dx -\int_\Gamma\Big( \meansc{\mathcal{G}\nabla w} \cdot\jumpsc{v}+\theta\meansc{\mathcal{G}\nabla v} \cdot\jumpsc{w}-\sigma\jumpsc{w}\cdot \jumpsc{v} \Big) \dS,
\]
and
\[
	B_{\rm a}(w,v):=\sum_{K\in\mathcal{T}}\bigg(-\int_{K}\mathbf{b} w\cdot \nabla v \dx +\int_{\partial_+K\cap\Gamma_b}(\mathbf{b}\cdot\mathbf{n})w^+v^+ \dS+\int_{\partial_-K\backslash\Gamma_b}(\mathbf{b}\cdot\mathbf{n})w^-\lfloor v\rfloor \dS\bigg).
\]
with $\partial_-K:=\{x\in\partial K: \mathbf{b}(x)\cdot \mathbf{n}(x)<0\}$ denoting the inflow part of the element boundary $\partial K$. For later use we note the identity
\begin{equation}
	\label{DG_adv_HJ_form}
	B_{\rm a}(w,v)= \sum_{K\in\mathcal{T}}\Big(\int_{K}\nabla\cdot (\mathbf{b} w)v\dx-\int_{\partial_-K\backslash\Gamma_b} \!\! (\mathbf{b} \cdot\mathbf{n})\lfloor w\rfloor v^+\dS
	-\int_{\partial_-K\cap\Gamma_b} \!\! (\mathbf{b} \cdot\mathbf{n}) w^+v^+\dS\Big),
\end{equation}
which can be shown by performing element-wise integration by parts and regrouping the face integrals appearing; we refer to \cite{newpaper} for details.

We define the respective IPDG norm $\ndg{w}:=\big(\ndg{w}_{\rm d}^2+\ndg{w}_{\rm a}^2\big)^{1/2}$, whereby
\[
	\begin{aligned}
		\ndg{w}_{\rm d}:= & \ \Big(
		\frac{1}{2}\|\sqrt{\mathcal{G}}\nabla_h w\|_{\Omega}^2
		+ \frac{1}{2}\|\sqrt{\sigma}\jumpsc{w}\|_{\Gamma}^2
		\Big)^{\frac{1}{2}}, \quad\text{and}\quad
		\ndg{w}_{\rm a}:= \Big(\frac{1}{2}\|\sqrt{|\mathbf{b}\cdot\mathbf{n}|}\jumpsc{w}\|_{\Gamma}^2\Big)^{\frac{1}{2}}.
	\end{aligned}
\]
Selecting $\sigma$ as in \eqref{eq:penalty_final}, and employing the trace inverse inequality \eqref{eq:trace_inverse_estimate}, we get \eqref{inverse_magic}, which, in turn, implies the coercivity and continuity of the bilinear form $B_{\rm d}$, viz.,
\begin{equation}\label{coer_d}
	\ndg{V}_{\rm d}^2\le B_{\rm d}(V,V),\qquad\text{and}\qquad 	 B_{\rm d}(V,W)\le 3\ndg{V}_{\rm d}\ndg{W}_{\rm d},
\end{equation}
for all $V,W\in V_h$. Moreover, we have the coercivity identity
\begin{equation}\label{coer_a}
	\ndg{V}_{\rm a}^2-\frac{1}{2}\int_\Omega(\nabla_h\cdot \mathbf{b}) V^2\dx= B_{\rm a}(V,V),
\end{equation}
for all $V\in V_h$, whose proof can be found, e.g., in \cite{newpaper}.

If we are to ensure stability in the $L_2$-norm, heuristically speaking, the maximum eigenvalue $\lambda_{\max}$ of $B$ is inversely proportional to the CFL constant.  To estimate the numerical range, which contains the eigenvalues, we form the Rayleigh quotient $\sup_{v\in V_h}\frac{B(V,V)}{\|V\|_{\Omega}^2}$. Then,
using \eqref{coer_d} and \eqref{coer_a} we have, respectively,
\begin{equation}\label{rayleigh}
	\ndg{V}^2-\frac{1}{2}\|\sqrt{|\nabla_h\cdot \mathbf{b}|} V\|_{L_2(\Omega)}^2
	\le B_{\rm d}(V,V)+B_{\rm a}(V,V)=B(V,V)
	\le 3\ndg{V}_{\rm d}^2+\ndg{V}_{\rm a}^2+\frac{1}{2}\|\sqrt{|\nabla_h\cdot \mathbf{b}|} V\|_{L_2(\Omega)}^2.
\end{equation}

Set $\beta_K :=\|\mathbf{b}\|_{L_\infty(K)}$ and $\beta'_K:= \frac{1}{2}\|\nabla\cdot\mathbf{b}\|_{L_\infty(K)}$, for brevity, as well as $G_K:=\max_{K'}\|\mathcal{G}\|_{L_\infty(K')}$, where the maximum is taken over all $K$ and all its face-neighbours.

We consider the trace-inverse inequality \eqref{eq:trace_inverse_estimate} and we set $C_{{\rm inv},K}=\max_{e\subset\partial K} C_{\rm inv}(e,K,p)$ for brevity, so that $\|V\|_e^2\le C_{{\rm inv},K}\|V\|_K^2$ for $v\in \mathbb{P}_K^p$. We also recall the classical inverse estimate
\[
	\|\nabla V\|_{K}^2\le C_{\rm inv,2}(K)\|V\|_{K}^2,\quad \text{ with }\quad C_{\rm inv,2}(K)=C_\nabla p^4\rho_K^{-2},
\] with $\rho_K$ the radius of the largest inscribed ball in $K$, for some generic global constant $C_\nabla>0$. Note that $\rho_K\sim {\rm diam}(K)$ for shape-regular mesh families. For later use we also set $C_{{\rm inv},2,K}:=\max_{K'}C_{\rm inv,2}(K')$ with the maximum taken over $K$ itself and all its face-neighbours.

Using the above, we have
\[
	\begin{aligned}
		\ndg{V}_{\rm d}^2\le & \ \frac{1}{2}\sum_{K\in\mathcal{T}}\|\sqrt{\mathcal{G}}\nabla_h V\|_{K}^2+\sum_{K\in\mathcal{T}}\max_{e\subset \partial K}\sigma_e\|V|_K\|_{e}^2 \\
		\le                  & \ \frac{1}{2}\sum_{K\in\mathcal{T}}\|\sqrt{\mathcal{G}}\nabla_h V\|_{K}^2+\sum_{e\subset\Gamma}\sum_{*\in\{+,-\}}\sigma|_e \|V|_{K^*}\|_{e}^2
		\le \sum_{K\in\mathcal{T}} \Big(\frac{1}{2}C_{\rm inv,2}(K) G_K
		+(d+1)C_{{\rm inv},K}\sigma_K\}\Big) \|V\|_{K}^2,
	\end{aligned}
\]
with $\sigma_K:=\max_{e\subset\partial K} \sigma_e$,
since each element has $d+1$ faces. The above gives
\begin{equation}\label{diff_bound}
	\ndg{V}_{\rm d}^2
	\le\ \sum_{K\in\mathcal{T}}\big(\frac{1}{2}C_{\rm inv, 2}(K)G_K+(d+1)C_{{\rm inv}, K}\sigma_K\big)\|V\|_{K}^2.
\end{equation}
(Note that $\sigma_K$ is, roughly speaking, proportional to $C_1G_KC_{{\rm inv}, K}$, with $C_1$ the penalty constant; cf., \eqref{eq:penalty_final}.)
Working analogously, we also have
\begin{equation}
	\label{adv_bound}
	\ndg{V}_{\rm a}^2\le \sum_{K\in\mathcal{T}}(d+1)\beta_K C_{{\rm inv}, K}\|V\|_{K}^2.
\end{equation}
Thus, \eqref{rayleigh} implies the following upper bound for the numerical range:
\begin{equation}\label{rayleigh_bound}
	\sup_{v\in V_h}\frac{B(V,V)}{\|V\|_{\Omega}^2}\le \max_{K\in\mathcal{T}}   \Big(\frac{3}{2}C_{\rm inv, 2}(K)G_K+(d+1)C_{{\rm inv}, K}(3\sigma_K+\beta_K)+\beta_K'\Big) =: \Lambda(\mathcal{G},\mathbf{b}, \mathcal{T},p) \equiv \Lambda.
\end{equation}

A typical choice for the CFL condition is to require $
	|\lambda_{\rm max}|\delta t\le C_{\rm CFL}
$,
for a sufficiently small positive constant $C_{\rm CFL}$, with $\lambda_{\max}$ denoting the largest eigenvalue in size of the discrete spatial operator. Thus, it is sufficient to consider a practical CFL condition of the form
$
	\Lambda \delta t\le C_{\rm CFL}$.
From the above, we can see that the size of the penalty parameter affects the maximum timestep $\delta t$ in an inversely proportional manner.

Alternatively, following \cite{burman_ern_fernandez},  the proof of stability for the classical strong stability-preserving third-order Runge-Kutta method therein relies on a bound of the form
\begin{equation}
	\label{BEF_claim}
	\|\mathcal{L}_hV\|_{\Omega}\le \tilde{\Lambda}\|V\|_{\Omega},
\end{equation}
for $\tilde{\Lambda}$, depending on the PDE coefficients, on the mesh and on $p$, whereby the discrete spatial operator (stiffness matrix) $\mathcal{L}_h:V_h\to V_h$, is defined by
\begin{equation}\label{riesz}
	\int_\Omega \mathcal{L}_h V W\dx = B(V,W), \qquad \text{for all } V,W\in V_h.
\end{equation}
(We note that \cite{burman_ern_fernandez} is concerned with purely transport problems, i.e., $\mathcal{G}=0$.) Once \eqref{BEF_claim} is established, the required CFL condition reads
$
	\tilde{\Lambda}\delta t\le C_{\rm CFL}$,
for a sufficiently small positive constant $C_{\rm CFL}$.
Thus, it suffices to seek an upper bound for $\tilde{\Lambda}$. To that end, \eqref{coer_d} gives
\begin{equation}\label{bound_L}
	\|\mathcal{L}_hV\|_{\Omega}^2 =B(V,\mathcal{L}_h V)\le |B_{\rm d}(V,\mathcal{L}_h V)|+|B_{\rm a}(V,\mathcal{L}_h V)|.
\end{equation}

For $B_{\rm a}(V,\mathcal{L}_h V)$, starting from \eqref{DG_adv_HJ_form}, we estimate each term as  follows:
\[
	\begin{aligned}
		|	B_{\rm a}(V,\mathcal{L}_h V)|
		\le & \
		\sum_{K\in\mathcal{T}_h}\Big(\big(\beta_K\|\nabla V\|_{K} + \beta_K' \|V\|_{K}\big) \|\mathcal{L}_h V\|_{K}
		+\beta_K\|V\|_{L_2(\partial_-K)}\|\mathcal{L}_hV\|_{L_2(\partial_-K)}\Big)                                                                                \\
		\le & \   \sum_{K\in\mathcal{T}_h}\Big(\big( \sqrt{C_{\rm inv,2}(K)}+(d+1)C_{{\rm inv}, K}\big)\beta_K+ \beta_K'\Big)\| V\|_{K}  \|\mathcal{L}_h V\|_{K},
	\end{aligned}
\]
upon employing inverse estimates, as above. Finally, Young's inequality of the form $\alpha\beta\le \alpha^2+\frac{1}{4}\beta^2$, yields
\begin{equation}\label{advection_CFL_one}
	\begin{aligned}
		|	B_{\rm a}(V,\mathcal{L}_h V)|
		\le & \   \sum_{K\in\mathcal{T}_h}\tilde{\Lambda}_{{\rm a},K}\| V\|_{K}^2 + \frac{1}{4} \|\mathcal{L}_h V\|_{\Omega}^2,
	\end{aligned}
\end{equation}
with
\[
	\tilde{\Lambda}_{{\rm a},K}:=\Big(\big( \sqrt{C_{\rm inv,2}(K)}+(d+1)C_{{\rm inv}, K}\big)\beta_K+ \beta_K'\Big)^2.
\]

For $B_{\rm d}(V,\mathcal{L}_h V)$, we begin by setting $\tilde{G}|_K:=\max_{K'}G_{K'}$, with the maximum taken over $K$ itself and all its face-neighbours. Working completely analogously to \eqref{inverse_magic}, we arrive at the bound
\[
	\begin{aligned}
		|B_{\rm d}(V,\mathcal{L}_h V)|
		\le & \sum_{K\in\mathcal{T}}\Big(G_K \|\nabla V\|_{K} \|\nabla\mathcal{L}_h V\|_{K} +\frac{\epsilon_1}{2}\tilde{G}_K\|\nabla V\|_{K}^2 +\frac{\epsilon_2\theta}{2}\tilde{G}_K\|\nabla \mathcal{L}_hV\|_{K}^2\Big) \\
		    & +\frac{1}{2}\sum_{e\subset\Gamma}\sigma_eC_1^{-1}\Big(\frac{1}{\epsilon_2}\|\jumpsc{V}\|_{e}^2+\frac{1}{\epsilon_1}\|\jumpsc{\mathcal{L}_h V}\|_{e}^2\Big),
	\end{aligned}
\]
employing Young's inequality $\alpha\beta\le \alpha^2/(2\epsilon_*)+\epsilon_*\beta^2/2$, for $\epsilon_*>0$, with $C_1$ as in the definition of $\sigma_e$ in \eqref{eq:penalty_final}.
Employing again inverse estimates on the last bound, we get
\[
	\begin{aligned}
		|B_{\rm d}(V,\mathcal{L}_h V)|
		\le & \sum_{K\in\mathcal{T}}\bigg(\tilde{G}_K C_{{\rm inv},2}(K)\Big(\| V\|_{K} \|\mathcal{L}_h V\|_{K} +\frac{\epsilon_1}{2}\|V\|_{K}^2 +\frac{\epsilon_2}{2}\theta\| \mathcal{L}_hV\|_{K}^2\Big) \\
		    & \qquad+\sigma_KC_1^{-1}C_{{\rm inv},K}\Big(\epsilon_2^{-1}\|V\|_{K}^2+\epsilon_1^{-1}\|\mathcal{L}_h V\|_{K}^2\Big)\bigg).
	\end{aligned}
\]
Employing Young's inequality, we have $\| V\|_{K} \|\mathcal{L}_h V\|_{K} \le \frac{\epsilon_3}{2}\| V\|_{K}^2+\frac{1}{2\epsilon_3}\|\mathcal{L}_h V\|_{K}^2$, for $\epsilon_3>0$, and factorising, we arrive at
\[
	\begin{aligned}
		|B_{\rm d}(V,\mathcal{L}_h V)|
		\le & \sum_{K\in\mathcal{T}}\bigg(\Big(\frac{1}{2}\tilde{G}_K C_{{\rm inv},2}(K)\big(\epsilon_1+\epsilon_3\big)+\frac{\sigma_KC_{{\rm inv},K}}{C_1\epsilon_2}\Big)\|V\|_{K}^2            \\
		    & \qquad+\Big(\frac{1}{2}\tilde{G}_K C_{{\rm inv},2}(K)\big(\epsilon_2\theta+\epsilon_3^{-1}\big)+\frac{\sigma_KC_{{\rm inv},K}}{C_1\epsilon_1}\Big)\| \mathcal{L}_hV\|_{K}^2\bigg).
	\end{aligned}
\]
We select $\epsilon_1= C_1\big(4\sigma_K C_{{\rm inv},K}\big)^{-1}$, $\epsilon_2=\big(2\tilde{G}_K C_{{\rm inv},2}(K)\big)^{-1} $, and $\epsilon_3=\epsilon_2^{-1}$, to deduce
\begin{equation}\label{diffusion_CFL_one}
	\begin{aligned}
		|B_{\rm d}(V,\mathcal{L}_h V)|
		\le & \sum_{K\in\mathcal{T}}\tilde{\Lambda}_{{\rm d},K}\|V\|_{K}^2 +\frac{1}{4}\| \mathcal{L}_hV\|_{\Omega}^2.
	\end{aligned}
\end{equation}
with
\[
	\tilde{\Lambda}_{{\rm d},K}:=\frac{C_1\tilde{G}_K C_{{\rm inv},2}(K)}{8\sigma_K C_{{\rm inv},K}}+\tilde{G}_K^2 C_{{\rm inv},2}^2(K)+\frac{2\tilde{G}_K\sigma_K}{C_1}C_{{\rm inv},K} C_{{\rm inv},2}(K).
\]

Applying \eqref{advection_CFL_one} and \eqref{diffusion_CFL_one} to \eqref{bound_L}, along with elementary calculations give \eqref{BEF_claim} with
\[
	\tilde{\Lambda}\equiv \tilde{\Lambda}(\mathcal{G},\mathbf{b}, \mathcal{T},p):=\max_{K\in\mathcal{T}}\sqrt{ 2(\tilde{\Lambda}_{{\rm a},K}+\tilde{\Lambda}_{{\rm d},K})}.
\]
Comparing the formulas for $\Lambda$ and $\tilde{\Lambda}$ we can see the similarities in their dependence upon the mesh-size and the polynomial degree. In particular, on shape-regular mesh families, we have $\Lambda\sim \tilde{\Lambda}$ with similarity constants independent of the PDE coefficients, the mesh, and the polynomial degree $p$.


\bibliographystyle{plain}

\bibliography{cas-refs.bib}


\end{document}